\documentclass[a4paper,12pt]{article}
\usepackage{amsmath}
\usepackage{amsthm}
\usepackage{amsfonts}
\usepackage{amssymb}
\usepackage{amstext}
\usepackage{amsopn}
\usepackage{graphicx}
\usepackage{mathrsfs}

\usepackage{color}
\usepackage{mathtools}

\newcommand{\al}{\alpha}		
\newcommand{\ze}{\zeta}		
\newcommand{\Gam}{\Gamma}	

\newcommand{\bC}{\mathbb C}		
\newcommand{\bN}{\mathbb N}		
\newcommand{\bR}{\mathbb R}		
\newcommand{\bZ}{\mathbb Z}		
\newcommand{\cA}{\mathcal A}		%
\newcommand{\cE}{\mathcal E}		
\newcommand{\cF}{\mathcal F}		
\newcommand{\cG}{\mathcal G}		
\newcommand{\cM}{\mathcal M}		
\newcommand{\cN}{\mathcal N}		
\newcommand{\sS}{\mathscr S}		%
\newcommand{\vC}{{\bf C}}		
\newcommand{\vL}{{\bf L}}		
\newcommand{\vT}{{\bf T}}		
\newcommand{\vW}{{\bf W}}	
\newcommand{\vZ}{{\bf Z}}		
\newcommand{\ve}{{\bf e}}		
\newcommand{\vr}{{\bf r}}		
\newcommand{\vz}{{\bf z}}		
\newcommand{\vze}{{\bf 0}}		
\newcommand{\val}{{\text{\boldmath{$\al$}}}}		
\newcommand{\dd}{{\rm d}}		
\newcommand{\ii}{{\rm i}}		
\newcommand{\re}{{\rm Re}}		
\newcommand{\im}{{\rm Im}}		
\newcommand{\Conv}{{\rm Conv}}		
\newcommand{\vol}{{\rm Vol}}		

\newcommand{\abs}[1]{\left|{#1}\right|}			
\newcommand{\rob}[1]{\left({#1}\right)}		
\newcommand{\cub}[1]{\left\{{#1}\right\}}		
\newcommand{\fr}{\frac}		
\newcommand{\tfr}{\tfrac}		
\newcommand{\nth}[1]{\fr{1}{#1}}			
\newcommand{\half}{\fr{1}{2}}			
\newcommand{\nonum}{\nonumber}	
\newcommand{\wt}{\widetilde}				
\newcommand{\ul}{\underline}				
\newcommand{\ol}{\overline}				
\newcommand{\dst}{\displaystyle}			
\newcommand{\tst}{\textstyle}				
\definecolor{r}{rgb}{0,0,0}					
\definecolor{g}{rgb}{0,0,0}					
\definecolor{b}{rgb}{0,0,0}					
\definecolor{p}{rgb}{0,0,0}					

\theoremstyle{definition}
	\newtheorem{dfn}{Definition}
	\newtheorem{eg}[dfn]{Example}
	\newtheorem*{rmk}{Remark}
\theoremstyle{plain}
	\newtheorem{lem}[dfn]{Lemma}
	\newtheorem{pro}[dfn]{Proposition}
	\newtheorem{thm}[dfn]{Theorem}

\begin{document}

\title{\bf Finiteness of Fixed Equilibrium Configurations of Point Vortices in the Plane with Background Flow}
\author{\textsc{Pak-Leong Cheung${}^{1}$\thanks{
{\color{b} Partially supported by a graduate studentship of HKU and the RGC grants HKU 706411P and HKU 703313P.}
\mbox{\hspace{4truecm}}} 
\ and Tuen Wai Ng${}^{1}$\thanks{
{\color{b}Partially supported by the RGC grant HKU 703313P.}
\mbox{\hspace{0.9truecm}}}
}}
 
\date{}
\maketitle

\centerline{\textbf{\today}}

\begin{figure}[b]
\rule[-2.5truemm]{5cm}{0.1truemm}\\[2mm]
{\footnotesize  
{\color{b}2000 {\it Mathematics Subject Classification: Primary} 70F10, 76B99, {\it Secondary} 13P15.}
\par {\it Key words and phrases.} 
{\color{b}point vortex, fixed equilibrium, polynomial system, BKK theory}
\par\noindent 1.
Department of Mathematics,
The University of Hong Kong,
Pokfulam, Hong Kong.
\smallskip
\par\noindent E-mail: {\color{b}\tt mathcpl@connect.hku.hk}, {\tt ntw@maths.hku.hk}}
\end{figure}

\begin{quotation}
\noindent{\textbf{Abstract}. For {\color{p}a dynamic} system consisting of $n$ point vortices in an ideal plane fluid with {\color{r}a steady, incompressible and} irrotational background flow, a more physically significant definition of {\color{p}a} fixed equilibrium configuration is suggested. Under this new definition, if the complex polynomial $w$ that determines the aforesaid background flow is non-constant, we have found an attainable generic upper bound \smash{$\fr{(m+n-1)!}{(m-1)!\,n_1!\cdots n_{i_0}!}$} for the number of fixed equilibrium configurations. Here, $m=\deg w$, $i_0$ is the number of species, and each $n_i$ is the number of vortices in a species. We transform the rational function system {\color{r}arisen from} fixed {\color{r}equilibria} into {\color{p}a polynomial} system, whose form is {\color{p}good} enough to apply the BKK theory (named after D. N. Bernshtein \cite{Ber}, A. G. Khovanskii \cite{Kh} and A. G. Kushnirenko \cite{Kus}) to show the finiteness of its number of solutions. Having this finiteness, the required bound follows from B\'ezout's theorem {\color{r}or the BKK root count by T. Y. Li and X.-S. Wang \cite{LiWa}}.}
\end{quotation}

\bigskip
\section{Introduction}
\label{sec1}

\bigskip
Every polynomial $w$ in one complex variable generates a {\color{r}{\it steady} (i.e. independent of time)} fluid flow in the complex plane $\bC$ by  the map $\ze\mapsto(\re\,w(\ze),-\im\,w(\ze))$, which is identified with $\ol{w(\ze)}$. The polynomial $w$ is called the {\it complex velocity} of this flow; and since $w$ satisfies the Cauchy--Riemann equations, the flow is {\color{r}{\it incompressible} and {\it irrotational}}. Assume that the fluid in $\bC$ is {\it ideal} (in the sense of \cite[Section 1.3, p.6]{Ac}). If $n\geq2$ (point) vortices $z_1(t),\dots,z_n(t)\in\bC$ with their respective {\it circulations}
\begin{equation}\label{eq4.1.1}
\Gam_1,\dots,\Gam_n\in\bR_*:=\bR\setminus\cub{0}
\end{equation}
(which are constants, due to Helmholtz's theorems (cf. \cite[Section 2.2]{MeAr})) are situated in the ideal fluid with the above flow $\ol w$ in the background, then the dynamics of these vortices will be governed by
\begin{equation}\label{eq4.1.2}
\fr{\dd z_j(t)}{\dd t}=-\nth{2\pi\ii}\sum_{\substack{k=1\\k\neq j}}^n\fr{\Gam_k}{\ \ol{z_j(t)}-\ol{z_k(t)}\ }+\ol{w(z_j(t))},\ j=1,\dots,n
\end{equation}
(cf. \cite[Equation (2)]{KadCa} (caution: typo on the left-hand side) and \cite[Equation (34)]{Cl}).

\bigskip
The case of no background flow is when $w\equiv0$. Since the early 1980s, there {\color{p}have been many studies on} all types of stationary configurations (in the sense of \cite[Definitions 0.2 \& 1.1.2]{ONe1}, see also \cite[Definitions 0.3 \& 1.1.3 \& Proposition 1.1.4]{ONe1}), as well as other aspects of the vortex dynamics \eqref{eq4.1.2}. So far, researchers have only found configurations with special patterns or for small $n$ (see \cite{CaKad}, \cite{Cl}, \cite{KadCa}, etc., for these results). Given the difficulty of determining (analytically or numerically) the stationary configurations, inquiry into {\color{p}their number} becomes a natural alternative. O'Neil \cite[Theorems 5.1.1, 5.2.1 \& 6.5.1]{ONe1} (cf. \cite[Propositions 1 \& 2]{ONe2}, and beware of different terminologies) gave such results for three types:
\begin{quote}
\it (Let $w\equiv0$.) For almost every choice of circulations \eqref{eq4.1.1} that satisfies
\begin{enumerate}
\item[$\bullet$]	$\sum_{i<j}\Gam_i\Gam_j=0$, there are exactly $(n-2)!$ equilibrium configurations. (This relation among the circulations \eqref{eq4.1.1} is necessary for the existence of such configurations.)
\item[$\bullet$]	$\sum_j\Gam_j=0$, there are exactly $(n-1)!$ rigidly translating configurations. (This relation among the circulations \eqref{eq4.1.1} is necessary for the existence of such configurations.)
\item[$\bullet$]	$\sum_{i<j}\Gam_i\Gam_j\neq0$ and $\sum_j\Gam_j\neq0$, there are no more than $n!/2$ collinear relative equilibrium configurations.
\end{enumerate}
\end{quote}
See also Hampton \cite{Ha} and Hampton and Moeckel \cite{HaMo2} for more results about {\color{p}the} number of configurations when $n=4,5$.

\bigskip
As for {\color{p}the case} when a background flow is present (i.e. $w\not\equiv0$), there appears to be no corresponding knowledge so far apart from the few cases in \cite[Sections III.B.2 and III.D]{KadCa} and \cite[Section 3.3]{Cl}. Our main result (Theorem \ref{thm4.6}) concerning the {\color{r}finiteness of the} number of {\it fixed equilibrium configurations} ({\color{r}to be defined in Definition \ref{dfn4.4}) will fill this deficiency. Also note that this terminology is synonymous with `equilibrium configurations' in \cite{ONe1} and `stationary equilibrium configurations' in \cite{ONe2}, but the present one seems more common in the literature.}

\bigskip
\section{{\color{r}Results}}
\label{sec2}

\bigskip
Before proceeding, two conventions are to be understood throughout this article: (i) `number of solutions' of any single polynomial or rational function equation, or any such system, counts multiplicity; (ii) `finitely many' includes `none'.

\bigskip
{\color{g}{\color{r}{\it Fixed equilibria} are the solutions of the rational function system \eqref{eq4.2.2} below in the $n$ unknowns $z_1,\dots,z_n$. Equivalently, they are the distinct solutions of the polynomial system $\sS$ obtained by clearing denominators in \eqref{eq4.2.2}. The inquiry into finiteness of the number of solutions of polynomial systems is reminiscent of a tool `reduced system test' (Lemma \ref{lem4.3}) in the {\it BKK theory} (as detailed in \cite[Section 3]{HaMo1}). However, this test does not work at least in some cases of this rather complicated polynomial system $\sS$. Therefore, we have found an alternative one \eqref{eq4.2.3} (Lemma \ref{lem4.1}) whose form is {\color{p}good} enough for applying the test to confirm the finiteness of the number of solutions of this new polynomial system \eqref{eq4.2.3} (Proposition \ref{pro4.2}).} These finitely many solutions $(z_1,\dots,z_n)$ are then reduced to {\color{r}\it fixed equilibrium configurations} by the notion of {\it equivalent solutions} (Definition \ref{dfn4.4}, which is different from those in \cite{HaMo2}, \cite{ONe1} and \cite{ONe2}). We will consider a natural situation (Definition \ref{dfn4.5}) where equivalent solutions arise, and arrive at the main result {\color{r}(Theorem \ref{thm4.6}). The rest of this section is devoted to expanding this paragraph.}}

\bigskip
For {\it fixed {\color{r}equilibria}}, set $z_j(t)\equiv z_j$ (so that \smash{$\fr{\dd z_j(t)}{\dd t}\equiv0$}) in the system \eqref{eq4.1.2} and hence we have
\begin{equation}\label{eq4.2.1}
-w(z_j)=\nth{2\pi\ii}\sum_{\substack{k=1\\k\neq j}}^n\fr{\Gam_k}{z_j-z_k},\ j=1,\dots,n.
\end{equation}
{\color{r}The case} $w\equiv c\in\bC_*:=\bC\setminus\cub{0}$ goes back to O'Neil's results in \cite{ONe1} and \cite{ONe2}, therefore we assume that the {\it degree $m:=\deg w$ of the background flow} is positive in what follows. Then, {\color{r}by complex conjugation and} an appropriate rescaling, it only suffices to consider the normalized system
\begin{equation}\label{eq4.2.2}
{z_j}^m+W(z_j)=\sum_{\substack{k=1\\k\neq j}}^n\fr{\Gam_k}{z_j-z_k}=:L_j(z_1,\dots,z_n),\ j=1,\dots,n,
\end{equation}
where $W$ is a polynomial of degree at most $m-1$ with coefficients determined by $w$. {\color{r}One} might clear denominators to obtain a polynomial system {\color{r}$\sS$} of $n$ equations, where each equation is of degree $m+n-1$ and is in the $n$ unknowns $z_1,\dots,z_n$. {\color{r}Polynomial systems could have infinitely or finitely many solutions. If $\sS$ falls into the latter case, then B\'ezout's theorem (\cite[Theorem 2.3.1]{Cox}) would provide
\begin{equation}\label{eq4.2.101}
(m+n-1)^n
\end{equation}
as an upper bound} for the number of solutions.

\bigskip
{\color{r}As far as finiteness of the number of solutions of $\sS$ is concerned, the following test in the BKK theory may be called upon (for details, the reader is referred to \cite[Section 3]{HaMo1}, and also \cite{HaMo2} for an application):}

\bigskip
{\color{g}\begin{lem}\label{lem4.3}
{\bf({\color{r}`Reduced system test'} for finiteness of {\color{p}the} number of solutions of {\color{p}a} polynomial system in {\color{p}$\pmb{\bC_*^n}$})}\quad
{\rm\cite[Propositions 2 \& 3]{HaMo1}}\quad Consider a system of $m$ polynomial equations in $n$ unknowns:
\begin{equation}\label{eq4.2.102}
P_k(z_1,\dots,z_n)=\sum_{\vr\,=(r_1,\dots,r_n)\in\cA_k}c_\vr{z_1}^{r_1}\cdots{z_n}^{r_n}=0,\quad c_\vr\in\bC_*,\quad k=1,\dots,m,
\end{equation}
where each finite subset $\cA_k$ of $(\bN\cup\cub{0})^n$ is called the {\it support} of $P_k$. For each $\val\in\bR^n$, the system
\begin{equation}\label{eq4.2.103}
P_{k,\val}(z_1,\dots,z_n):=\sum_{\vr\in\cA_k,\ \val\cdot\vr=\min_{\vr'\in\cA_k}\val\cdot\vr'}c_\vr{z_1}^{r_1}\cdots{z_n}^{r_n}=0,\quad k=1,\dots,m,
\end{equation}
is called the reduced system (of \eqref{eq4.2.102}) determined by $\val$. Suppose that there exists an $\val_0\in\bZ^n$ such that every reduced system \eqref{eq4.2.103} with $\val_0\cdot\val\leq0$ has no solution in $\bC_*^n$. Then, the original system \eqref{eq4.2.102} has only finitely many solutions in $\bC_*^n$.
\end{lem}

\bigskip
\begin{rmk}\quad
The seemingly weird condition `$\val\cdot\vr=\min_{\vr'\in\cA_k}\val\cdot\vr'$' in \eqref{eq4.2.103} actually admits a beautiful geometric interpretation in terms of {\it supporting hyperplane}. We refer the reader to the paragraph preceding \cite[Proposition 3]{HaMo1}. This interpretation will be used when proving Proposition \ref{pro4.2} in Section \ref{sec4}.
\end{rmk}}

\bigskip
{\color{r}Despite the availability of such {\color{p}a} handy finiteness test, the structure of $\sS$ is still too complicated for the test to conclude anything even in some cases with small $n$ (the number of vortices) and $m$ (the degree of the background flow). More precisely, there are reduced systems which do have solutions in $\bC_*^n$ but no choice of $\val_0$ could avoid all of these reduced systems.} This difficulty has motivated us to transform \eqref{eq4.2.2} into the following better polynomial system:

\bigskip
\begin{lem}\label{lem4.1}
{\bf(An equivalent polynomial system)}\quad
The system \eqref{eq4.2.2} is equivalent to
\begin{equation}\label{eq4.2.3}
\left\{\begin{alignedat}{4}
&\sum_j\Gam_j{z_j}^m		&&+\sum_j\Gam_jW(z_j)			&&=0\\
&\sum_j\Gam_j{z_j}^{m+1}		&&+\sum_j\Gam_jz_jW(z_j)		&&=\sum_{i<j}\Gam_{i,j}\\
&\sum_j\Gam_j{z_j}^{m+2}		&&+\sum_j\Gam_j{z_j}^2W(z_j)		&&=\sum_j\Gam_j\Gam^jz_j\\
&						&&							&&\ \,\vdots\\
&\sum_j\Gam_j{z_j}^{m+k-1}	&&+\sum_j\Gam_j{z_j}^{k-1}W(z_j)	&&=\sum_j\Gam_j\Gam^j{z_j}^{k-2}+\sum_{\substack{r+s=k-2\\r,s\neq0,\ i<j}}\Gam_{i,j}{z_i}^r{z_j}^s\\
&						&&							&&\ \,\vdots\\
&\sum_j\Gam_j{z_j}^{m+n-1}	&&+\sum_j\Gam_j{z_j}^{n-1}W(z_j)	&&=\sum_j\Gam_j\Gam^j{z_j}^{n-2}+\sum_{\substack{r+s=n-2\\r,s\neq0,\ i<j}}\Gam_{i,j}{z_i}^r{z_j}^s
\end{alignedat}\right.
\end{equation}
with constraint
\begin{equation}\label{eq4.2.4}
z_i\neq z_j,\ \ i\neq j,
\end{equation}
where $\Gam_{i,j}:=\Gam_i\Gam_j$ and $\Gam^j:={\sum_{i\neq j}}\Gam_i$.
\end{lem}

\bigskip
\begin{rmk}\quad
\begin{enumerate}
\item[(i)]	The $n$ equations {\color{r}of this new system \eqref{eq4.2.3}} are of degrees $m$, $m+1$, ... and $m+n-1$ respectively, and the leftmost sums are the only sources of these degrees.
\item[(ii)]	The right-hand side of the {\color{r}$k$-th} ($k=2,\dots,n$) equation {\color{r}of \eqref{eq4.2.3}} is either a homogeneous polynomial of degree $k-2$ or, in the degenerate case, identically zero.
\item[(iii)]	Each equation {\color{r}of \eqref{eq4.2.3}} is invariant under any permutation of the $n$ circulation-vortex pairs $\cub{(\Gam_j,z_j):j=1,\dots,n}$, while in the original system \eqref{eq4.2.2}, there will be a permutation of the $n$ equations.
\item[(iv)]	{\it Admissible} solutions of \eqref{eq4.2.3} are those that satisfy constraint \eqref{eq4.2.4}. This terminology will not appear until Section \ref{sec5}.
\item[(v)]	{\color{r}Lemma \ref{lem4.1} will be proved in Section \ref{sec3} by transforming \eqref{eq4.2.2} into \eqref{eq4.2.3} via an explicit matrix \eqref{matrix} which is invertible under \eqref{eq4.2.4}. This new system \eqref{eq4.2.3} {\color{b}can also be obtained by} replacing $2\pi\ii\bar v_j$ by ${z_j}^m+W(z_j)$ in O'Neil's \cite[Equation (4.1)]{ONe2}. Our and O'Neil's methods are {\color{b}different}. In particular, the involvement of the Vandermonde determinant in \eqref{eq4.3.5} in our method seems surprising.}
\end{enumerate}
\end{rmk}

\bigskip\noindent
{\color{g}As the reader will see in Section \ref{sec4}, the advantage of Lemma \ref{lem4.1}'s transformation is that \eqref{eq4.2.3} is in a special form that facilitates using {\color{r}Lemma \ref{lem4.3}.}}

\bigskip
The following proposition shows that the condition
\begin{equation}\label{eq4.2.5}
\sum_{j\in I}\Gam_j\neq0\text{ \ for all \ }I\subset\cub{1,\dots,n}\quad\text{and}\quad\sum_{i<j}\Gam_i\Gam_j\neq0
\end{equation}
guarantees that \eqref{eq4.2.3} alone (i.e. without the constraint \eqref{eq4.2.4}) already has finitely many solutions:

\bigskip
\begin{pro}\label{pro4.2}
{\bf(Upper bound for {\color{p}the} number of solutions of {\color{p}the} equivalent polynomial system)}\quad
Assume that $m\geq1$. Then, for any choice of the circulations \eqref{eq4.1.1} that satisfies \eqref{eq4.2.5}, the system \eqref{eq4.2.3} has at most \smash{$\fr{(m+n-1)!}{(m-1)!}$} solutions, and so does \eqref{eq4.2.2}. This bound can be attained.
\end{pro}

\bigskip
\begin{rmk}\quad
There are cases where \eqref{eq4.2.2} has infinitely many solutions, such as in \cite[Section III.B \& III.C]{KadCa} and \cite[Section 3.2]{Cl}. But these existing cases either have no background flow or a background flow of degree $0$ (i.e. $m<1$), thus are not covered by Proposition \ref{pro4.2}.
\end{rmk}

\bigskip\noindent
{\color{r}In Proposition \ref{pro4.2}, the} finiteness of the number of solutions of \eqref{eq4.2.3} under assumption \eqref{eq4.2.5} will be shown in Section \ref{sec4}. Then by Remark (i) {\color{r}following Lemma \ref{lem4.1}}, B\'ezout's theorem provides the required upper bound \smash{$\fr{(m+n-1)!}{(m-1)!}$} for the number of solutions which is better than the bound \eqref{eq4.2.101}. As we shall see in the proof of Lemma \ref{lem4.1} in Section \ref{sec3}, every solution of \eqref{eq4.2.2} also satisfies \eqref{eq4.2.3}, so \eqref{eq4.2.2} inherits the upper bound for its number of solutions. Finally, Section \ref{sec5} will provide examples where this bound is attained.

\bigskip
Next, we suggest a perhaps new definition of {\color{p}a} {\it fixed equilibrium configuration} in order to state the main result Theorem \ref{thm4.6} of this article. Fix any polynomial $w$ of degree $m\geq1$ which provides the background flow in $\bC$ as in Section \ref{sec1}. Let $\vz=(z_1,\dots,z_n)$ be a solution of \eqref{eq4.2.1} (necessarily, $z_i\neq z_j$ for $i\neq j$). Recall that these vortices in the background flow $w$ generate the flow
$$\ol V_\vz(\ze):=-\nth{2\pi\ii}\sum_j\fr{\Gam_j}{\ \ol\ze-\ol{z_j}\ }+\ol{w(\ze)}$$
or the complex velocity
$$V_\vz(\ze):=\ol{\ol V_\vz(\ze)}=\nth{2\pi\ii}\sum_j\fr{\Gam_j}{\ze-z_j}+w(\ze)$$
on $\bC_\vz:=\bC\setminus\cub{z_1,\dots,z_n}$. Note that $V_\vz$ is a rational function on $\bC$ with simple poles at $z_1,\dots,z_n$.

\bigskip
\begin{dfn}\label{dfn4.4}
{\bf(Equivalent solutions. Fixed equilibrium configuration)}\quad\newline
Two solutions $\vz=(z_1,\dots,z_n)$ and $\vz'=(z'_1,\dots,z'_n)$ of \eqref{eq4.2.1} are said to be {\it equivalent} (denoted by $\vz\sim\vz'$) if they generate two non-reflectively similar flows with $w$, in a plane-geometrical sense. More precisely, $\vz\sim\vz'$ if there exists $(a,b)\in\bC_*\times\bC$ such that \smash{$\ol a\ol V_{\vz'}(a\ze+b)=\ol V_\vz(\ze)$} for $\ze\in\bC_\vz$, or, {\color{p}equivalently (in terms of complex velocity),}
\begin{equation}\label{eq4.2.6}
aV_{\vz'}(a\ze+b)=V_\vz(\ze)\text{ \ for \ }\ze\in\bC_\vz.
\end{equation}
Then, each equivalence class $[\vz]$ is called a {\it fixed equilibrium configuration}.
\end{dfn}

\bigskip\noindent
In O'Neil's definition \cite[Definitions 0.4 \& 1.1.5]{ONe1}, two solutions $\vz$ and $\vz'$ are regarded as equivalent if there exists $(a,b)\in\bC_*\times\bC$ such that $z_j'=az_j+b$ for all $j$, so only the shapes of vortex sets are involved. Our definition of equivalent solutions of \eqref{eq4.2.1} has an extra physical significance: Besides the shapes of vortex sets, their effects on the rest of the plane are also considered. Here we illustrate the difference between our and O'Neil's definitions:

\bigskip
\begin{eg}\label{eg4.101}
\quad
Consider the two-vortex case $n=2$, both with circulation $\Gam_1=\Gam_2=1$, and with background flow \smash{$w(\ze)=-\fr{\ze^2+1}{2\pi\ii}$} (thus $m=2$). By Lemma \ref{lem4.1}, we are to solve
$$\left\{\begin{alignedat}{4}
&{z_1}^2+{z_2}^2+2		&&=0\\
&{z_1}^3+{z_2}^3+z_1+z_2	&&=1
\end{alignedat}\right.$$
with constraint $z_1\neq z_2$, and the solutions are
\begin{equation}\label{eq4.2.104}
\begin{alignedat}{2}
(z_1,z_2)\approx&\ (-0.250\pm1.349\,\ii,0.487\pm0.693\,\ii),\\
&\ (0.487\pm0.693\,\ii,-0.250\pm1.349\,\ii),\\
&\ (-0.237\pm1.028\,\ii,-0.237\mp1.028\,\ii).
\end{alignedat}
\end{equation}
Any pair of two-point sets in the plane differ by translation, rotation and/or scaling, so all the six solutions in \eqref{eq4.2.104} are equivalent in O'Neil's sense, thereby constituting exactly one fixed equilibrium configuration. But by Definition \ref{dfn4.4}, they constitute three fixed equilibrium configurations because they generate three flows with $w$ that are not non-reflectively similar as shown in Figure \ref{fig1}. The streamlines in Figure \ref{fig1} are actually formed by superimpositions in Figure \ref{fig2}. (Figures \ref{fig1} and \ref{fig2} are generated by Mathematica.)
\begin{figure}[h]
\centering
\begin{tabular}{@{}cc@{}}
\includegraphics[width=.43\textwidth]{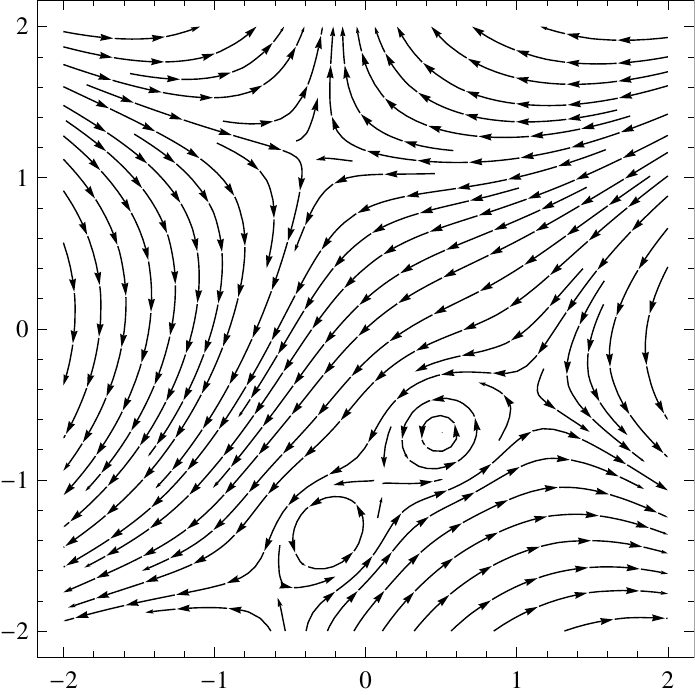}&\includegraphics[width=.43\textwidth]{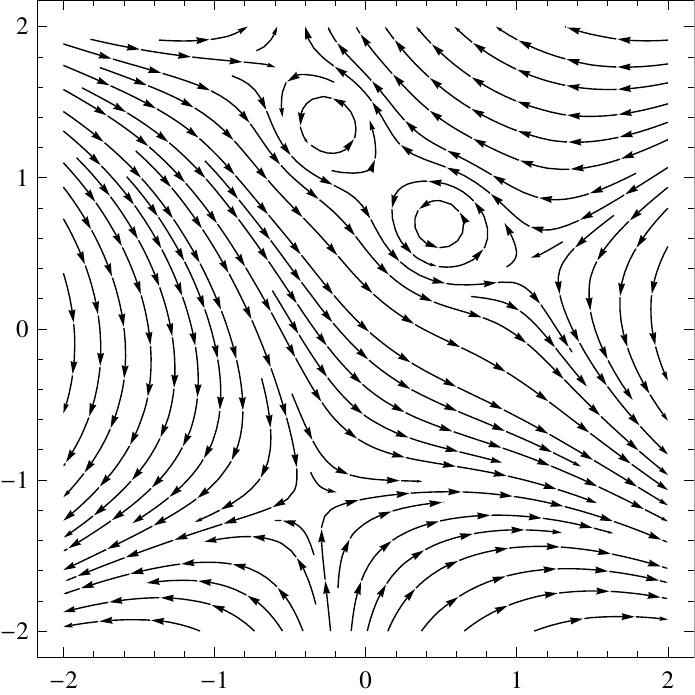}\\
\scriptsize(a) $(z_1,z_2)\approx(-0.250-1.349\,\ii,0.487-0.693\,\ii)$&\scriptsize(b) $(z_1,z_2)\approx(-0.250+1.349\,\ii,0.487+0.693\,\ii)$
\end{tabular}\vspace{.3cm}\\
\begin{tabular}{@{}c@{}}
\includegraphics[width=.43\textwidth]{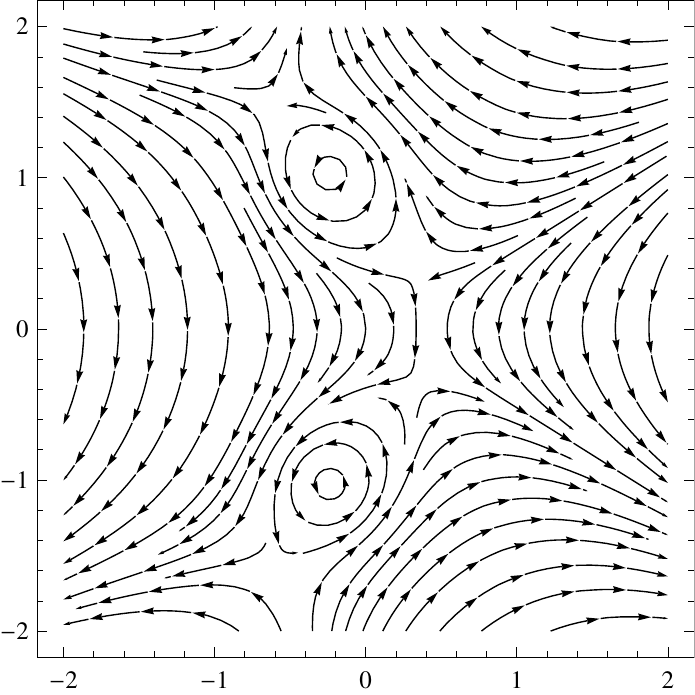}\\
\scriptsize(c) $(z_1,z_2)\approx(-0.237-1.028\,\ii,-0.237+1.028\,\ii)$
\end{tabular}
\caption{The streamlines of the flows \smash{$\ol V_{(z_1,z_2)}$} generated by the vortex sets \eqref{eq4.2.104} with the background flow \smash{$w(\ze)=-\fr{\ze^2+1}{2\pi\ii}$} in Example \ref{eg4.101}.\label{fig1}}
\end{figure}
\begin{figure}[h]
\centering
\begin{tabular}{@{}cc@{}}
\includegraphics[width=.43\textwidth]{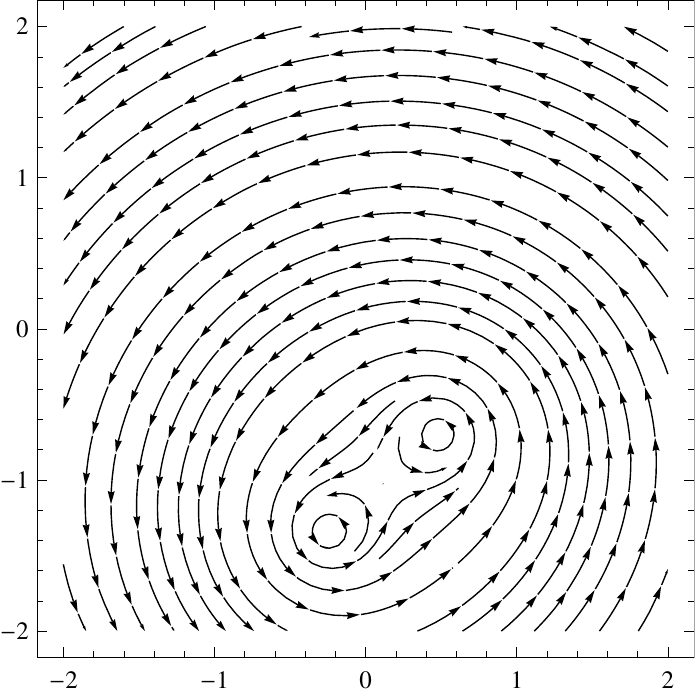}&\includegraphics[width=.43\textwidth]{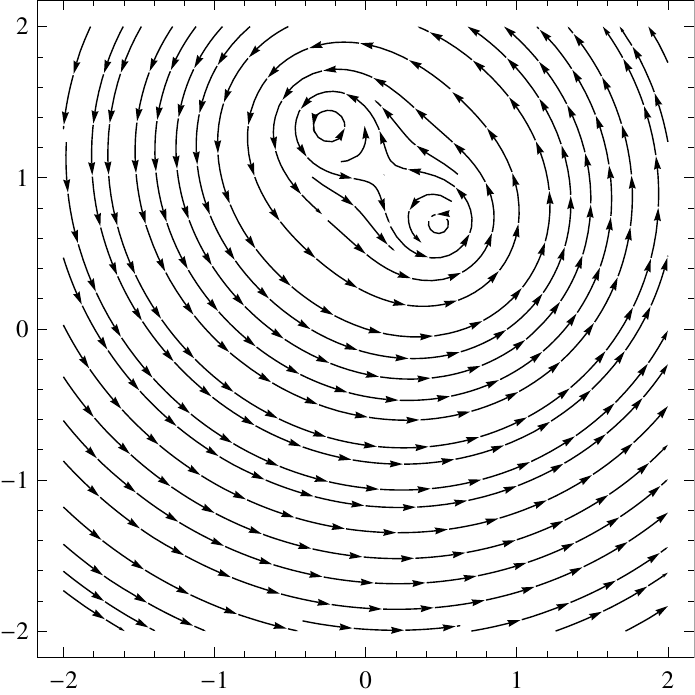}\\
\scriptsize(a) $(z_1,z_2)\approx(-0.250-1.349\,\ii,0.487-0.693\,\ii)$&\scriptsize(b) $(z_1,z_2)\approx(-0.250+1.349\,\ii,0.487+0.693\,\ii)$\vspace{.3cm}\\
\includegraphics[width=.43\textwidth]{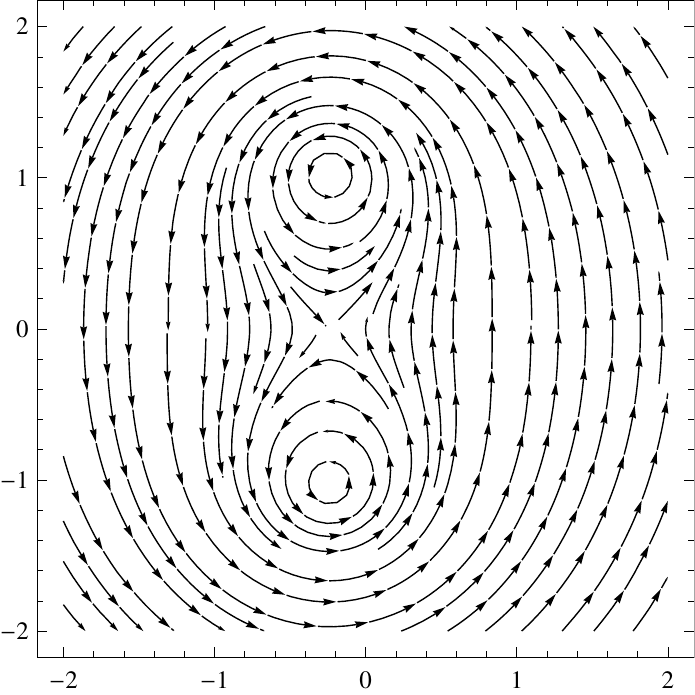}&\includegraphics[width=.43\textwidth]{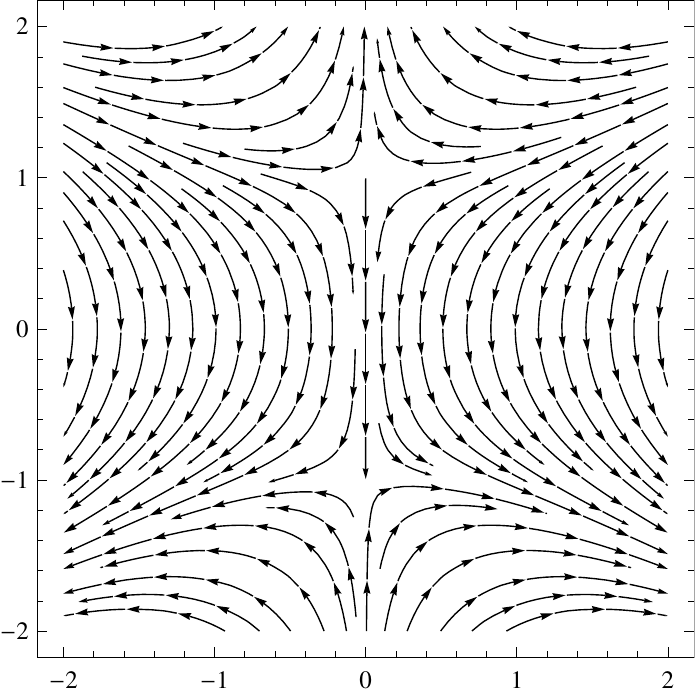}\\
\scriptsize(c) $(z_1,z_2)\approx(-0.237-1.028\,\ii,-0.237+1.028\,\ii)$&\scriptsize(d) background flow
\end{tabular}
\caption{The streamlines of the flows generated by the vortex sets \eqref{eq4.2.104} without the background flow \smash{$w(\ze)=-\fr{\ze^2+1}{2\pi\ii}$} in Example \ref{eg4.101}. The background flow (d) superimposes with (a)--(c) respectively to form {\color{p}Figures} \ref{fig1}(a)--(c).\label{fig2}}
\end{figure}
\end{eg}

\bigskip
Our definition introduces comparison \eqref{eq4.2.6} of rational functions to distinguish between different configurations. We are not going to elaborate this definition to the full extent, but just to observe a natural situation where equivalent solutions, in our sense, of \eqref{eq4.2.1} arise:

\bigskip
\begin{dfn}\label{dfn4.5}
{\bf(Species)}\quad
Partition the circulations \eqref{eq4.1.1} into their $i_0$ distinct values:
\begin{equation}\label{eq4.2.7}
\Gam_{j_{i,1}}=\cdots=\Gam_{j_{i,n_i}}=:[\Gam_i],\ \ i=1,\dots,i_0,\ \ [\Gam_i]\neq[\Gam_{i'}],\ \ i\neq i',\ \ \tst\sum_in_i=n.
\end{equation}
For each solution $(z_1,\dots,z_n)$ of \eqref{eq4.2.1}, the set of all the vortices
$$[z_i]:=\cub{z_{j_{i,1}},\dots,z_{j_{i,n_i}}}$$
which possess the common circulation $[\Gam_i]$ is called a {\it species} (cf. \cite{KadCa}).
\end{dfn}

\bigskip\noindent
Note that the concept of species was not involved in defining equivalent solutions in \cite{ONe1}. Here, it acts as follow: If $\Gam_i=\Gam_j$ for some $i<j$, then {\color{p}with} every solution $\vz=(z_1,\dots,z_i,\dots,z_j,\dots,z_n)$ of \eqref{eq4.2.1} there associates another solution $\vz'=(z_1,\dots,z_j,\dots,z_i,\dots,z_n)$, and then $\vz\sim\vz'$ because \eqref{eq4.2.6} holds with $(a,b)=(1,0)$. Consequently, a standard combinatorial argument bridges Proposition \ref{pro4.2} to our main theorem:

\bigskip
\begin{thm}\label{thm4.6}
{\bf(Upper bound for {\color{p}the} number of fixed equilibrium configurations {\color{p}in a} background flow)}\quad
Assume that $m\geq1$. Then, for any choice of the circulations \eqref{eq4.1.1} that satisfies \eqref{eq4.2.5}, there are at most \smash{$\fr{(m+n-1)!}{(m-1)!\,n_1!\cdots n_{i_0}!}$} fixed equilibrium configurations, where $n_i$ are the sizes of the species as in \eqref{eq4.2.7}. This bound can be attained.
\end{thm}

\bigskip\noindent
Section \ref{sec5} will provide examples where this bound is and is not attained. It will also suggest another factor in reducing the solutions of \eqref{eq4.2.1}, \eqref{eq4.2.2} or \eqref{eq4.2.3} to fixed equilibrium configurations via Definition \ref{dfn4.4}.

\bigskip
\section{Proof of Lemma \ref{lem4.1}}
\label{sec3}

\bigskip
Write \eqref{eq4.2.2} in vector form
\begin{equation}\label{eq4.3.1}
\vZ_n+\vW_n=\vL_n,
\end{equation}
where
\begin{equation}\label{eq4.3.2}
\vZ_n:=\begin{pmatrix}
{z_i}^m
\end{pmatrix}_{i=1,\dots,n},\ \ \vW_n:=\begin{pmatrix}
W(z_i)
\end{pmatrix}_{i=1,\dots,n}\text{ \ and \ }\vL_n:=\begin{pmatrix}
L_i(z_1,\dots,z_n)
\end{pmatrix}_{i=1,\dots,n}.
\end{equation}
Left-multiplying each side by the square matrix
\begin{equation}\label{matrix}
\vT_n:=\begin{pmatrix}
\Gam_j{z_j}^{i-1}
\end{pmatrix}_{i,j=1,\dots,n},
\end{equation}
the left-hand side simply becomes
\begin{align}
\vT_n(\vZ_n+\vW_n)&=\begin{pmatrix}
\begin{alignedat}{4}
&\Gam_1			&&\quad\Gam_2		&&\quad\cdots	&&\quad\Gam_n\\
&\Gam_1z_1		&&\quad\Gam_2z_2		&&\quad\cdots	&&\quad\Gam_nz_n\\
&\ \vdots			&&\quad\ \vdots			&&\quad\ddots	&&\quad\ \vdots\\
&\Gam_1{z_1}^{n-1}	&&\quad\Gam_2{z_2}^{n-1}	&&\quad\cdots	&&\quad\Gam_n{z_n}^{n-1}
\end{alignedat}
\end{pmatrix}\rob{\begin{pmatrix}
{z_1}^m\\
\vdots\\
{z_n}^m
\end{pmatrix}+\begin{pmatrix}
W(z_1)\\
\vdots\\
W(z_n)
\end{pmatrix}}\nonum\\
&=\begin{pmatrix}
\begin{alignedat}{1}
&\sum_j^{\phantom{n}}\Gam_j{z_j}^m\\
&\sum_j^{\phantom{n}}\Gam_j{z_j}^{m+1}\\
&\qquad\quad\vdots\\
&\sum_j^{\phantom{n}}\Gam_j{z_j}^{m+n-1}
\end{alignedat}
\end{pmatrix}+\begin{pmatrix}
\begin{alignedat}{1}
&\sum_j^{\phantom{n}}\Gam_jW(z_j)\\
&\sum_j^{\phantom{n}}\Gam_jz_jW(z_j)\\
&\qquad\quad\vdots\\
&\sum_j^{\phantom{n}}\Gam_j{z_j}^{n-1}W(z_j)
\end{alignedat}
\end{pmatrix}.\label{eq4.3.4}
\end{align}
On the right-hand side, the first entry of $\vT_n\vL_n$ is
$$\begin{pmatrix}
\Gam_1&\cdots&\Gam_n
\end{pmatrix}\begin{pmatrix}
\dst\sum_{i\neq 1}^{\phantom{n}}\fr{\Gam_i}{z_1-z_i}\\
\dst\sum_{i\neq 2}^{\phantom{n}}\fr{\Gam_i}{z_2-z_i}\\
\vdots\\
\dst\sum_{i\neq n}^{\phantom{n}}\fr{\Gam_i}{z_n-z_i}
\end{pmatrix}=\sum_j\Gam_j\sum_{i\neq j}\fr{\Gam_i}{z_j-z_i}=\sum_{i<j}\rob{\fr{\Gam_j\Gam_i}{z_j-z_i}+\fr{\Gam_i\Gam_j}{z_i-z_j}}=0$$
by cancellations. The second entry is
\begin{align*}
\begin{pmatrix}
\Gam_1z_1&\cdots&\Gam_nz_n
\end{pmatrix}\begin{pmatrix}
\dst\sum_{i\neq 1}^{\phantom{n}}\fr{\Gam_i}{z_1-z_i}\\
\dst\sum_{i\neq 2}^{\phantom{n}}\fr{\Gam_i}{z_2-z_i}\\
\vdots\\
\dst\sum_{i\neq n}^{\phantom{n}}\fr{\Gam_i}{z_n-z_i}
\end{pmatrix}&=\sum_j\Gam_jz_j\sum_{i\neq j}\fr{\Gam_i}{z_j-z_i}\\
&=\sum_{i<j}\rob{\fr{\Gam_j\Gam_iz_j}{z_j-z_i}+\fr{\Gam_i\Gam_jz_i}{z_i-z_j}}=\sum_{i<j}\Gam_i\Gam_j.
\end{align*}
The {\color{r}$k$-th} entry, $k=3,\dots,n$, is
\begin{align*}
\begin{pmatrix}
\Gam_1{z_1}^{k-1}&\cdots&\Gam_n{z_n}^{k-1}
\end{pmatrix}&\begin{pmatrix}
\dst\sum_{i\neq 1}^{\phantom{n}}\fr{\Gam_i}{z_1-z_i}\\
\dst\sum_{i\neq 2}^{\phantom{n}}\fr{\Gam_i}{z_2-z_i}\\
\vdots\\
\dst\sum_{i\neq n}^{\phantom{n}}\fr{\Gam_i}{z_n-z_i}
\end{pmatrix}=\sum_j\Gam_j{z_j}^{k-1}\sum_{i\neq j}\fr{\Gam_i}{z_j-z_i}\\
&\qquad=\sum_{i<j}\rob{\fr{\Gam_j\Gam_i{z_j}^{k-1}}{z_j-z_i}+\fr{\Gam_i\Gam_j{z_i}^{k-1}}{z_i-z_j}}\\
&\qquad=\sum_{i<j}\Gam_i\Gam_j({z_i}^{k-2}+{z_i}^{k-3}z_j+\cdots+z_i{z_j}^{k-3}+{z_j}^{k-2}).
\end{align*}
We have left-multiplied each side of \eqref{eq4.3.1} by $\vT_n$ to obtain \eqref{eq4.2.3}, thus the claimed equivalence between \eqref{eq4.2.2} and `\eqref{eq4.2.3} and \eqref{eq4.2.4}' would follow if $\vT_n$ is invertible. Indeed,
\begin{equation}\label{eq4.3.5}
\det\vT_n=\abs{\begin{alignedat}{4}
&\Gam_1			&&\quad\Gam_2		&&\quad\cdots		&&\quad\Gam_n\\
&\Gam_1z_1		&&\quad\Gam_2z_2		&&\quad\cdots		&&\quad\Gam_nz_n\\
&\quad\vdots		&&\quad\quad\vdots		&&\quad\ddots	&&\quad\quad\vdots\\
&\Gam_1{z_1}^{n-1}	&&\quad\Gam_2{z_2}^{n-1}	&&\quad\cdots		&&\quad\Gam_n{z_n}^{n-1}
\end{alignedat}}=\prod_j\Gam_j\cdot\prod_{i<j}(z_j-z_i)
\end{equation}
is just \smash{$\prod_j\Gam_j$} times of the Vandermonde determinant \smash{$\prod_{i<j}(z_j-z_i)$}, where \eqref{eq4.1.1} and \eqref{eq4.2.4} guarantees that their product is non-zero.

\bigskip
\section{Proof of Proposition \ref{pro4.2}}
\label{sec4}

\bigskip
\subsection{A Little More General System} 
\label{sec4.1}

\bigskip
Consider the following more general system than \eqref{eq4.2.3}:
\begin{equation}\label{eq4.4.1}
\left\{\begin{alignedat}{5}
&P_1(z_1,\dots,z_n)	&&:=\sum_j\Gam_j{z_j}^m		&&+\sum_{j,r}A^{1,r}_j{z_j}^r												&&=0\\
&P_2(z_1,\dots,z_n)	&&:=\sum_j\Gam_j{z_j}^{m+1}	&&+\sum_{j,r}A^{2,r}_j{z_j}^{r+1}+C_0											&&=0\\
&P_3(z_1,\dots,z_n)	&&:=\sum_j\Gam_j{z_j}^{m+2}	&&+\sum_{j,r}A^{3,r}_j{z_j}^{r+2}+\sum_jC^3_jz_j									&&=0\\
&				&&\ \ \,\vdots\\
&P_k(z_1,\dots,z_n)	&&:=\sum_j\Gam_j{z_j}^{m+k-1}	&&+\sum_{j,r}A^{k,r}_j{z_j}^{r+k-1}\\
&				&&						&&+\sum_jC^k_j{z_j}^{k-2}+\sum_{\substack{r+s=k-2\\r,s\neq0,\ i<j}}C^{k,r,s}_{i,j}{z_i}^r{z_j}^s	&&=0\\
&				&&\ \ \,\vdots\\
&P_n(z_1,\dots,z_n)	&&:=\sum_j\Gam_j{z_j}^{m+n-1}	&&+\sum_{j,r}A^{n,r}_j{z_j}^{r+n-1}\\
&				&&						&&+\sum_jC^n_j{z_j}^{n-2}+\sum_{\substack{r+s=n-2\\r,s\neq0,\ i<j}}C^{n,r,s}_{i,j}{z_i}^r{z_j}^s	&&=0,
\end{alignedat}\right.
\end{equation}
with coefficients $\Gam_j,C_0\in\bC_*$ and $A_{r,j},C^k_j,C^{k,r,s}_{i,j}\in\bC$, where \smash{$\sum_r=\sum_{r=0}^{m-1}$}. We will prove the following finiteness result for \eqref{eq4.4.1} via Sections \ref{sec4.2} and \ref{sec4.3}:

\bigskip
\begin{pro}\label{pro4.7}
{\bf(Finiteness of {\color{p}the} number of solutions)}\quad
The system \eqref{eq4.4.1} has {\color{p}only} finitely many solutions if \smash{$\sum_{j\in I}\Gam_j\neq0$} for all $I\subset\cub{1,\dots,n}$ and $C_0\neq0$.
\end{pro}

\bigskip\noindent
The system \eqref{eq4.4.1} has B\'ezout bound \smash{$\fr{(m+n-1)!}{(m-1)!}$}. This coincides with the generally finer Li and Wang's BKK root count \cite{LiWa}, due to the somewhat special structure of \eqref{eq4.4.1} (to be seen in \eqref{eq4.4.3}):

\bigskip
\begin{pro}\label{pro4.8}
{\bf(BKK root count in {\color{p}$\pmb{\bC^n}$})}\quad
Assume that $\Gam_j,C_0\in\bC_*$. Let $\cM(\cN_1,\dots,\cN_n)$ denote the mixed volume of the Newton polytopes $\cN_k$ of $\cA_k\cup\cub{\vze}$, $k=1,\dots,n$, where each $\cA_k$ is the support of the polynomial $P_k$ in \eqref{eq4.4.1}. Then,
$$\cM(\cN_1,\dots,\cN_n)=\fr{(m+n-1)!}{(m-1)!}.$$
\end{pro}

\bigskip
\begin{proof}[{\bf Proof of Proposition \ref{pro4.8}.}]\quad
The supports of $P_k$ are respectively
\begin{equation}\label{eq4.4.2}
\begin{alignedat}{2}
\cA_1	&=m\cE_n		\cup\bigcup_{j,r}S^{1,r}_j\\
\cA_2	&=(m+1)\cE_n	\cup\bigcup_{j,r}S^{2,r}_j	\cup\cub{\vze}\\
\cA_3	&=(m+2)\cE_n	\cup\bigcup_{j,r}S^{3,r}_j	\cup\bigcup_jS^3_j\\
		&\ \,\vdots\\
\cA_k	&=(m+k-1)\cE_n	\cup\bigcup_{j,r}S^{k,r}_j	\cup\bigcup_jS^k_j	\cup\bigcup_{\substack{r+s=k-2\\r,s\neq0,\ i<j}}S^{k,r,s}_{i,j}\\
		&\ \,\vdots\\
\cA_n	&=(m+n-1)\cE_n	\cup\bigcup_{j,r}S^{n,r}_j	\cup\bigcup_jS^n_j	\cup\bigcup_{\substack{r+s=n-2\\r,s\neq0,\ i<j}}S^{n,r,s}_{i,j},
\end{alignedat}
\end{equation}
where $\cE_n:=\cub{\ve_1,\dots,\ve_n}$ is the standard basis of $\bR^n$, \smash{$\bigcup_r=\bigcup_{r=0}^{m-1}$},
$$S^{k,r}_j:=\left\{\begin{matrix}
\cub{(r+k-1)\ve_j}	&\text{if}	&A^{k,r}_j\neq0\\
\emptyset			&\text{if}	&A^{k,r}_j=0
\end{matrix}\right.,\quad S^k_j:=\left\{\begin{matrix}
\cub{(k-2)\ve_j}	&\text{if}	&C^k_j\neq0\\
\emptyset		&\text{if}	&C^k_j=0
\end{matrix}\right.$$
$$\text{and}\quad S^{k,r,s}_{i,j}:=\left\{\begin{matrix}
\cub{r\ve_i+s\ve_j}	&\text{if}	&C^{k,r,s}_{i,j}\neq0\\
\emptyset			&\text{if}	&C^{k,r,s}_{i,j}=0
\end{matrix}\right..$$
No matter what {\color{p}value $S^{k,r}_j$}, $S^k_j$ and $S^{k,r,s}_{i,j}$ take in \eqref{eq4.4.2}, the Newton polytopes $\cN_k$ of $\cA_k\cup\cub{\vze}$ are
\begin{equation}\label{eq4.4.3}
\begin{alignedat}{2}
\cN_1	&=\Conv(m\cE_n\cup\cub{\vze})			=m\Delta_n\\
\cN_2	&=\Conv((m+1)\cE_n\cup\cub{\vze})		=(m+1)\Delta_n\\
\cN_3	&=\Conv((m+2)\cE_n\cup\cub{\vze})		=(m+2)\Delta_n\\
		&\ \,\vdots\\
\cN_k	&=\Conv\rob{(m+k-1)\cE_n\cup\cub{\vze}}	=(m+k-1)\Delta_n\\
		&\ \,\vdots\\
\cN_n	&=\Conv\rob{(m+n-1)\cE_n\cup\cub{\vze}}	=(m+n-1)\Delta_n,
\end{alignedat}
\end{equation}
where
\begin{equation}\label{eq4.4.4}
\Delta_n:=\Conv(\cE_n\cup\cub{\vze})
\end{equation}
is the unit simplex in $\bR^n$. Such {\color{p}a} simplification is due to that $(r+k-1)\ve_j$ and $r\ve_i+s\ve_j$ ($k=3,\dots,n$, $r+s=k-2$ and $i<j$) actually lie in $\Conv((m+k-1)\cE_n\cup\cub{\vze})$. By the multilinearity \cite[Theorem 7.4.12.b, p.338]{CoxLitOS2} of mixed volume and \cite[Exercise 7.3.b, p.306, \& Exercise 7.7.b, p.338]{CoxLitOS2}, the required mixed volume is
\begin{align*}
\cM(\cN_1,\dots,\cN_n)&=\cM(m\Delta_n,(m+1)\Delta_n,\dots,(m+n-1)\Delta_n)\\
&=m\cdot(m+1)\cdot\cdots\cdot(m+n-1)\cdot\cM(\Delta_n,\Delta_n,\dots,\Delta_n)\\
&=\fr{(m+n-1)!}{(m-1)!}\cdot n!\,\vol_n(\Delta_n)=\fr{(m+n-1)!}{(m-1)!}.\qedhere
\end{align*}
\end{proof}

\bigskip\noindent
Equation \eqref{eq4.4.3} has just revealed that the elements of the supports $\cA_k$ of $P_k$ corresponding to the terms with coefficients $A^\cdot_\cdot$ and $C^\cdot_\cdot$ can simply be ignored in the formation of the Newton polytopes $\cN_k$. This will be useful in the upcoming sections.

\bigskip
\subsection{The Initial Case}
\label{sec4.2}

\bigskip
To prove Proposition \ref{pro4.7}, we shall apply strong induction on the number of vortices $n$; and it contains Proposition \ref{pro4.2} as a case, and then leads to Theorem \ref{thm4.6} as analyzed in Section \ref{sec2}. Lemma \ref{lem4.3} will enter to test the finiteness of the number of solutions of polynomial systems in $\bC_*^n$. Taking the geometric interpretation of reduced system in Lemma \ref{lem4.3} as detailed in the paragraph preceding \cite[Proposition 3]{HaMo1} for granted, we will need the following notation for brevity of the upcoming discussions: let $\ell^\cN_\val$ denote the supporting hyperplane of a polytope $\cN$ with inward normal vector $\val$.

\bigskip
Here we start with the initial case $n=2$:
\begin{equation}\label{eq4.4.5}
\left\{\begin{matrix}
\begin{alignedat}{5}
&P_1(z_1,z_2)	&&=\Gam_1{z_1}^m+\Gam_2{z_2}^m+A^{1,m-1}_1{z_1}^{m-1}+A^{1,m-1}_2{z_2}^{m-1}\\
&			&&\qquad\qquad\qquad\qquad\qquad\qquad\qquad\quad+\cdots+(A^{1,0}_1+A^{1,0}_2)					&&=0\\
&P_2(z_1,z_2)	&&=\Gam_1{z_1}^{m+1}+\Gam_2{z_2}^{m+1}+A^{2,m-1}_1{z_1}^m+A^{2,m-1}_2{z_2}^m+\cdots+C_0	&&=0
\end{alignedat}\\
\Gam_1,\ \Gam_2,\ \Gam_1+\Gam_2,\ C_0\neq0
\end{matrix}\right.
\end{equation}

\bigskip\noindent
\ul{Case 1:\quad$(z_1,z_2)\in\bC_*^2$}.\quad Lemma \ref{lem4.3} will show that \eqref{eq4.4.5} has only finitely many solutions in $\bC_*^2$. To this end, consider, by \eqref{eq4.4.2}, the Newton polytopes
\begin{alignat*}{2}
\cN_1	&=\Conv(\cub{m\ve_1,m\ve_2}\cup\cdots)\quad\text{and}\\
\cN_2	&=\Conv(\cub{(m+1)\ve_1,(m+1)\ve_2,\vze}\cup\cdots)
\end{alignat*}
of the supports $\cA_1$ and $\cA_2$ of $P_1$ and $P_2$, and their Minkowski sum
\begin{equation}\label{eq4.4.6}
\cN=\cN_1+\cN_2=\Conv(\cub{a_2\ve_1,a_2\ve_2}\cup\cdots),
\end{equation}
where and hereafter,
\begin{equation}\label{eq4.4.7}
a_n:=m+(m+1)+\cdots+(m+n-1)=nm+\tst\binom{n}{2}.
\end{equation}
Note that {\color{p}this} `$\,\cdots$' in \eqref{eq4.4.6} does not alter the fact that
\begin{equation}\label{eq4.4.8}
\cN\subset a_2\Delta_2,
\end{equation}
where $\Delta_2$ is as in \eqref{eq4.4.4}. Now, consider the reduced systems of \eqref{eq4.4.5} determined by all the $\val$ or $\ell^\cN_\val$ satisfying $(\ve_1+\ve_2)\cdot\val\leq0$. Because of \eqref{eq4.4.8} and such {\color{p}a} choice of $\val$, each $\ell^\cN_\val$ actually supports $\cN$ at a face of the facet $\Conv(\cub{a_2\ve_1,a_2\ve_2})=a_2\Delta_2$:
\begin{enumerate}
\item[$\bullet$]	\ul{Case I:\quad$\ell^\cN_\val$ supports $\cN$ at the 0-face}
$$a_2\ve_j,\ \ j=1,2,$$
then $\ell^{\cN_1}_\val$ and $\ell^{\cN_2}_\val$ support $\cN_1$ and $\cN_2$ at the 0-faces
$$m\ve_j\quad\text{and}\quad(m+1)\ve_j$$
respectively, giving the reduced system
$$\left\{\begin{matrix}
\begin{alignedat}{4}
&\Gam_j{z_j}^m	&&=0\\
&\Gam_j{z_j}^{m+1}	&&=0
\end{alignedat}\\
\Gam_j\neq0
\end{matrix}\right.$$
which has no solution in $\bC_*^2$.
\item[$\bullet$]	\ul{Case II:\quad$\ell^\cN_\val$ supports $\cN$ at the 1-face}
$$\Conv(\cub{a_2\ve_1,a_2\ve_2}),$$
then $\ell^{\cN_1}_\val$ and $\ell^{\cN_2}_\val$ support $\cN_1$ and $\cN_2$ at the 1-faces
$$\Conv(\cub{m\ve_1,m\ve_2})\quad\text{and}\quad\Conv(\cub{(m+1)\ve_1,(m+1)\ve_2})$$
respectively, giving the reduced system
\begin{equation}\label{eq4.4.9}
\left\{\begin{matrix}
\begin{alignedat}{4}
&\wt P_1(z_1,z_2)	&&:=\Gam_1{z_1}^m+\Gam_2{z_2}^m			&&=0\\
&\wt P_2(z_1,z_2)	&&:=\Gam_1{z_1}^{m+1}+\Gam_2{z_2}^{m+1}	&&=0
\end{alignedat}\\
\Gam_1,\ \Gam_2,\ \Gam_1+\Gam_2\neq0.
\end{matrix}\right.
\end{equation}
\begin{enumerate}
\item[(i)]	If \ul{$z_1\neq z_2$}, then $\det\vT_2\neq0$ by \eqref{eq4.3.5}, so that $\vT_2$ is invertible. Left-multiplying both sides of \eqref{eq4.4.9} by ${\vT_2}^{-1}$, it is transformed to
$$\vZ_2=\vze$$
(see \eqref{eq4.3.4}), where $\vZ_2$ is as in \eqref{eq4.3.2}. But this contradicts that $z_1\neq z_2$.
\item[(ii)]	If \ul{$z_1=z_2$}, then the sub-system consisting of the first equation
$$\left\{\begin{matrix}
\begin{alignedat}{2}
&(\Gam_1+\Gam_2){z_1}^m&&=0
\end{alignedat}\\
\Gam_1+\Gam_2\neq0
\end{matrix}\right.$$
already has no solution in $\bC_*$.
\end{enumerate}
Therefore, in any case, \eqref{eq4.4.9} has no solution in $\bC_*^2$.
\end{enumerate}
What we have established so far is that for every $\val$ or $\ell^\cN_\val$ with $(\ve_1+\ve_2)\cdot\val\leq0$, the corresponding reduced system has no solution in $\bC_*^2$. Hence, it follows from Lemma \ref{lem4.3} that \eqref{eq4.4.5} has only finitely many solutions in $\bC_*^2$.

\bigskip\noindent
\ul{Case 2:\quad$(z_1,z_2)\in\bC_*\times\cub{0}$}.\quad The system \eqref{eq4.4.5} degrades to 
\begin{equation}\label{eq4.4.10}
\left\{\begin{matrix}
\begin{alignedat}{4}
&P_1(z_1,0)	&&=\Gam_1{z_1}^m+A^{1,m-1}_1{z_1}^{m-1}+\cdots+(A^{1,0}_1+A^{1,0}_2)	&&=0\\
&P_2(z_1,0)	&&=\Gam_1{z_1}^{m+1}+A^{2,m-1}_1{z_1}^m+\cdots+C_0				&&=0\\
\end{alignedat}\\
\Gam_1,\ C_0\neq0,
\end{matrix}\right.
\end{equation}
where the first equation already has only finitely many (at most $m$) solutions, thus so does \eqref{eq4.4.10}. Similar for the case where $(z_1,z_2)\in\cub{0}\times\bC_*$.

\bigskip\noindent
\ul{Case 3:\quad$(z_1,z_2)\in\cub{0}^2$}.\quad The system \eqref{eq4.4.5} degrades to
$$\left\{\begin{matrix}
\begin{alignedat}{6}
&P_1(0,0)	&&=A^{1,0}_1+A^{1,0}_2	&&=0\\
&P_2(0,0)	&&=C_0				&&=0\\
\end{alignedat}\\
C_0\neq0
\end{matrix}\right.$$
which is simply inconsistent.

\bigskip
Combining all the above three cases, \eqref{eq4.4.5} has only finitely many solutions in $\bC^2$, and Proposition \ref{pro4.7} with $n=2$ is proved.

\bigskip
\subsection{Strong Induction}
\label{sec4.3}

\bigskip
As the reader may have noticed or will see soon, everything actually lies in the non-solvability of the following special reduced system in $\bC_*^k$ under the assumption of Proposition \ref{pro4.7}:

\bigskip
\begin{lem}\label{lem4.9}
{\bf(Special reduced system)}\quad
The system
\begin{equation}\label{eq4.4.11}
\left\{\begin{matrix}
\begin{alignedat}{4}
&\wt P_1(z_1,\dots,z_k)	&&:=\Gam_1{z_1}^m+\cdots+\Gam_k{z_k}^m			&&=0\\
&\wt P_2(z_1,\dots,z_k)	&&:=\Gam_1{z_1}^{m+1}+\cdots+\Gam_k{z_k}^{m+1}		&&=0\\
&					&&\ \ \vdots\\
&\wt P_k(z_1,\dots,z_k)	&&:=\Gam_1{z_1}^{m+k-1}+\cdots+\Gam_k{z_k}^{m+k-1}	&&=0
\end{alignedat}\\
\sum_{j\in I}\Gam_j\neq0,\ \ I\subset\cub{1,\dots,k}
\end{matrix}\right.
\end{equation}
has no solution in $\bC_*^k$.
\end{lem}

\bigskip
\begin{proof}[{\bf Proof of Lemma \ref{lem4.9}}.]\quad
The case $k=2$ is just \eqref{eq4.4.9}. Assume that the lemma holds when $k=K$, then consider \eqref{eq4.4.11} with $k=K+1$:
\begin{equation}\label{eq4.4.12}
\left\{\begin{matrix}
\begin{alignedat}{4}
&\wt P_1(z_1,\dots,z_K,z_{K+1})	&&=\Gam_1{z_1}^m+\cdots+\Gam_K{z_K}^m+\Gam_{K+1}{z_{K+1}}^m			&&=0\\
&\wt P_2(z_1,\dots,z_K,z_{K+1})	&&=\Gam_1{z_1}^{m+1}+\cdots+\Gam_K{z_K}^{m+1}+\Gam_{K+1}{z_{K+1}}^{m+1}	&&=0\\
&							&&\ \ \vdots\\
&\wt P_K(z_1,\dots,z_K,z_{K+1})	&&=\Gam_1{z_1}^{m+K-1}+\cdots+\Gam_K{z_K}^{m+K-1}\\
&							&&\qquad\qquad\qquad\qquad\quad+\Gam_{K+1}{z_{K+1}}^{m+K-1}				&&=0\\
&\wt P_{K+1}(z_1,\dots,z_K,z_{K+1})	&&=\Gam_1{z_1}^{m+K}+\cdots+\Gam_K{z_K}^{m+K}\\
&							&&\qquad\qquad\qquad\qquad\quad+\Gam_{K+1}{z_{K+1}}^{m+K}	&&=0
\end{alignedat}\\
\sum_{j\in I}\Gam_j\neq0,\ \ I\subset\cub{1,\dots,K+1}
\end{matrix}\right.
\end{equation}
\begin{enumerate}
\item[(i)]	If \ul{all $z_j$ are distinct}, then $\det\vT_n\neq0$ by \eqref{eq4.3.5}, so that $\vT_n$ is invertible. Left-multiplying both sides of \eqref{eq4.4.12} by ${\vT_n}^{-1}$, it is transformed to
$$\vZ_n=\vze$$
(see \eqref{eq4.3.4}), where $\vZ_n$ is as in \eqref{eq4.3.2}. But this contradicts that all $z_j$ are distinct.
\item[(ii)]	If \ul{some $z_j$ are equal}, say, $z_K=z_{K+1}$, then the sub-system consisting of the first $K$ equations
$$\left\{\begin{matrix}
\begin{alignedat}{4}
&\wt P_1(z_1,\dots,z_K,z_K)	&&=\Gam_1{z_1}^m+\cdots+\Gam_{K-1}{z_{K-1}}^m+(\Gam_K+\Gam_{K+1}){z_K}^m			&&=0\\
&\wt P_2(z_1,\dots,z_K,z_K)	&&=\Gam_1{z_1}^{m+1}+\cdots+\Gam_{K-1}{z_{K-1}}^{m+1}\\
&							&&\qquad\qquad\qquad\qquad+(\Gam_K+\Gam_{K+1}){z_K}^{m+1}	&&=0\\
&							&&\ \ \vdots\\
&\wt P_K(z_1,\dots,z_K,z_K)	&&=\Gam_1{z_1}^{m+K-1}+\cdots+\Gam_{K-1}{z_{K-1}}^{m+K-1}\\
&							&&\qquad\qquad\qquad\qquad\quad+(\Gam_K+\Gam_{K+1}){z_K}^{m+K-1}						&&=0
\end{alignedat}\\
\sum_{j\in I}\Gam_j\neq0,\ \ I\subset\cub{1,\dots,K+1}
\end{matrix}\right.$$
already has no solution in $\bC_*^K$ due to the induction hypothesis (with $\Gam_K$ replaced by $\Gam_K+\Gam_{K+1}$), thus so does \eqref{eq4.4.12} in $\bC_*^{K+1}$. Similar for the other cases.
\end{enumerate}
Therefore, in any case, \eqref{eq4.4.12} has no solution in $\bC_*^{K+1}$, and the lemma follows from induction.
\end{proof}

\bigskip
Finally, we are in a position to prove Proposition \ref{pro4.7}. Recall that the case $n=2$ has already been proved above. Now, assume that it holds for $n=2,\dots,N-1$, then consider the case $n=N$:
\begin{equation}\label{eq4.4.13}
\left\{\begin{matrix}
\begin{alignedat}{5}
&P_1(z_1,\dots,z_N)	&&=\sum_j\Gam_j{z_j}^m		&&+\sum_{j,r}A^{1,r}_j{z_j}^r													&&=0\\
&P_2(z_1,\dots,z_N)	&&=\sum_j\Gam_j{z_j}^{m+1}	&&+\sum_{j,r}A^{2,r}_j{z_j}^{r+1}+C_0												&&=0\\
&P_k(z_1,\dots,z_N)	&&=\sum_j\Gam_j{z_j}^{m+k-1}	&&+\sum_{j,r}A^{k,r}_j{z_j}^{r+k-1}\\
&				&&						&&+\sum_jC^k_j{z_j}^{k-2}+\sum_{\substack{r+s=k-2\\r,s\neq0,\ i<j}}C^{k,r,s}_{i,j}{z_i}^r{z_j}^s		&&=0,
\end{alignedat}\\
k=3,\dots,N,\ \sum_{j\in I}\Gam_j\neq0,\ \ I\subset\cub{1,\dots,N},\ C_0\neq0
\end{matrix}\right.
\end{equation}

\bigskip\noindent
\ul{Case 1:\quad$(z_1,\dots,z_N)\in\bC_*^N$}.\quad Lemma \ref{lem4.3} will show that \eqref{eq4.4.13} has only finitely many solutions in $\bC_*^N$. To this end, consider, by \eqref{eq4.4.2}, the Newton polytopes
\begin{alignat*}{2}
&\cN_1	&&=\Conv(\cub{m\ve_1,\dots,m\ve_N}\cup\cdots)\\
&\cN_2	&&=\Conv(\cub{(m+1)\ve_1,\dots,(m+1)\ve_N,\vze}\cup\cdots)\\
&\cN_3	&&=\Conv(\cub{(m+2)\ve_1,\dots,(m+2)\ve_N}\cup\cdots)\\
&		&&\vdots\\
&\cN_N	&&=\Conv(\cub{(m+N-1)\ve_1,\dots,(m+N-1)\ve_N}\cup\cdots)
\end{alignat*}
of the supports $\cA_1,\dots,\cA_N$ of $P_1,\dots,P_N$, and their Minkowski sum of
$$\cN=\cN_1+\cdots+\cN_N=\Conv\rob{\cub{a_N\ve_1,\dots,a_N\ve_N}\cup\cdots},$$
where $a_N$ is as in \eqref{eq4.4.7}. Note that {\color{p}this} `$\,\cdots$' does not alter the fact that
$$\cN\subset a_N\Delta_N.$$
Now, consider the reduced systems of \eqref{eq4.4.13} determined by all the $\val$ or $\ell^\cN_\val$ satisfying $(\ve_1+\cdots+\ve_N)\cdot\val\leq0$. Because of \eqref{eq4.4.8} and such {\color{p}a} choice of $\val$, each $\ell^\cN_\val$ actually supports $\cN$ at a face of the facet $\Conv(\cub{a_N\ve_1,\dots,a_N\ve_N})=a_N\Delta_N$:
\begin{enumerate}
\item[$\bullet$]	\ul{Case I:\quad$\ell^\cN_\val$ supports $\cN$ at the 0-face}
$$a_N\ve_j,\ \ j=1,\dots,N,$$
then $\ell^{\cN_1}_\val,\dots,\ell^{\cN_N}_\val$ support $\cN_1,\dots,\cN_N$ at the 0-faces
$$m\ve_j,\dots,(m+N-1)\ve_j$$
respectively, giving the reduced system
$$\left\{\begin{matrix}
\begin{alignedat}{4}
&\Gam_j{z_j}^m 		&&=0\\
&\Gam_j{z_j}^{m+1}		&&=0\\
&					&&\vdots\\
&\Gam_j{z_j}^{m+N-1}	&&=0
\end{alignedat}\\
\Gam_j\neq0
\end{matrix}\right.$$
which has no solution in $\bC_*^N$.
\item[$\bullet$]	\ul{Case II:\quad $\ell^\cN_\val$ supports $\cN$ at the $(k-1)$-face ($k=2,\dots,N$)}
$$\Conv(\cub{a_N\ve_{j_1},\dots,a_N\ve_{j_k}}),\ \ j_1<\cdots<j_k,$$
then $\ell^{\cN_1}_\val,\dots,\ell^{\cN_N}_\val$ support $\cN_1,\dots,\cN_N$ at the $(k-1)$-faces
$$\Conv(\cub{m\ve_{j_1},\dots,m\ve_{j_k}}),\dots,\Conv(\cub{(m+N-1)\ve_{j_1},\dots,(m+N-1)\ve_{j_k}})$$
respectively, giving the reduced system
$$\left\{\begin{matrix}
\begin{alignedat}{4}
&\wt P_1(z_1,\dots,z_N)	&&=\Gam_{j_1}{z_{j_1}}^m+\cdots+\Gam_{j_k}{z_{j_k}}^m		&&=0\\
&\wt P_2(z_1,\dots,z_N)	&&=\Gam_{j_1}{z_{j_1}}^{m+1}+\cdots+\Gam_{j_k}{z_{j_k}}^{m+1}	&&=0\\
&					&&\ \ \vdots\\
&\wt P_k(z_1,\dots,z_N)	&&=\Gam_{j_1}{z_{j_1}}^{m+k-1}+\cdots+\Gam_{j_k}{z_{j_k}}^{m+k-1}	&&=0\\
&					&&\ \ \vdots\\
&\wt P_N(z_1,\dots,z_N)	&&=\Gam_{j_1}{z_{j_1}}^{m+N-1}+\cdots+\Gam_{j_k}{z_{j_k}}^{m+N-1}	&&=0,
\end{alignedat}\\
\sum_{j\in I}\Gam_j\neq0,\ \ I\subset\cub{j_1,\dots,j_k}
\end{matrix}\right.$$
where the sub-system consisting of the first $k$ equations already has no solution in $\bC_*^k$ by Lemma \ref{lem4.9}.
\end{enumerate}
What we have established so far is that for every $\val$ or $\ell^\cN_\val$ with $(\ve_1+\cdots+\ve_N)\cdot\val\leq0$, the corresponding reduced system has no solution in $\bC_*^N$. Hence, it follows from Lemma \ref{lem4.3} that \eqref{eq4.4.13} has only finitely many solutions in $\bC_*^N$.

\bigskip\noindent
\ul{Case 2:\quad $(z_1,\dots,z_N)\in\bC_*^k\times\cub{0}^{N-k}$ ($k=2,\dots,N-1$)}.\quad The system \eqref{eq4.4.13} degrades to
\begin{equation}\label{eq4.4.14}
\left\{\begin{matrix}
\begin{alignedat}{5}
&P_1(z_1,\dots,z_k,0,\cdots,0)	&&=\sum_{j=1}^k\Gam_j{z_j}^m+\sum_{j=1}^k\sum_rA^{1,r}_j{z_j}^r													&&=0\\
&P_2(z_1,\dots,z_k,0,\cdots,0)	&&=\sum_{j=1}^k\Gam_j{z_j}^{m+1}+\sum_{j=1}^k\sum_rA^{2,r}_j{z_j}^{r+1}+C_0												&&=0\\
&P_k(z_1,\dots,z_k,0,\cdots,0)	&&=\sum_{j=1}^k\Gam_j{z_j}^{m+k-1}+\sum_{j=1}^k\sum_rA^{k,r}_j{z_j}^{r+k-1}\\
&						&&\qquad+\sum_{j=1}^kC^k_j{z_j}^{k-2}+\sum_{\substack{r+s=k-2\\r,s\neq0,\ i<j}}C^{k,r,s}_{i,j}{z_i}^r{z_j}^s		&&=0,
\end{alignedat}\\
k=3,\dots,N,\ \sum_{j\in I}\Gam_j\neq0,\ \ I\subset\cub{1,\dots,k},\ C_0\neq0
\end{matrix}\right.
\end{equation}
where the sub-system consisting of the first $k$ equations already has only finitely many solutions in $\bC^k$ by induction hypothesis, thus so does \eqref{eq4.4.14}. Similar for the other cases where in $(z_1,\dots,z_N)$ exactly $N-k$ coordinates equal $0$.

\bigskip\noindent
\ul{Case 3:\quad $(z_1,\dots,z_N)\in\bC_*\times\cub{0}^{N-1}$}.\quad The system \eqref{eq4.4.13} degrades to
\begin{equation}\label{eq4.4.15}
\left\{\begin{matrix}
\begin{alignedat}{5}
&P_1(z_1,0,\cdots,0)	&&:=\Gam_1{z_1}^m		&&+\sum_rA^{1,r}_1{z_1}^r						&&=0\\
&P_2(z_1,0,\cdots,0)	&&:=\Gam_1{z_1}^{m+1}	&&+\sum_rA^{2,r}_1{z_1}^{r+1}+C_0				&&=0\\
&P_k(z_1,0,\cdots,0)	&&:=\Gam_1{z_1}^{m+k-1}	&&+\sum_rA^{k,r}_1{z_1}^{r+k-1}+C^k_1{z_1}^{k-2}	&&=0,
\end{alignedat}\\
k=3,\dots,N,\ \Gam_1\neq0
\end{matrix}\right.
\end{equation}
where the first equation already has only finitely many (at most $m$) solutions, thus so does \eqref{eq4.4.15}. Similar for the other cases where in $(z_1,\dots,z_N)$ exactly $N-1$ coordinates equal $0$.

\bigskip\noindent
\ul{Case 4:\quad $(z_1,\dots,z_N)\in\cub{0}^N$}.\quad The system \eqref{eq4.4.13} degrades to
$$\left\{\begin{matrix}
\begin{alignedat}{5}
&P_1(0,0,\cdots,0)	&&:=\sum_jA^{1,0}_j	&&=0\\
&P_2(0,0,\cdots,0)	&&:=C_0			&&=0\\
&P_k(0,0,\cdots,0)	&&:=0			&&=0,
\end{alignedat}\\
k=3,\dots,N,\ C_0\neq0
\end{matrix}\right.$$
is simply inconsistent.

\bigskip
Combining all the above four cases, \eqref{eq4.4.13} has only finitely many solutions in $\bC^N$, and Proposition \ref{pro4.7} follows from strong induction.

\bigskip
\section{Background Flow of Degree One}
\label{sec5}

\bigskip
By a background flow {\it of degree one}, we mean, without any loss of generality, \smash{$w(\ze)=-\fr{\ze+c}{2\pi\ii}$} for some constant $c\in\bC$ in \eqref{eq4.2.1} or $m=1$ and $W\equiv c$ in \eqref{eq4.2.2} or \eqref{eq4.2.3} (the case with $c=0$ is commonly called a {\it quadrupole background flow}). For whichever two distinct solutions $\vz=(z_1,\dots,z_n)$ and $\vz'=(z'_1,\dots,z'_n)$ of \eqref{eq4.2.1} with such $w$ to be equivalent, by \eqref{eq4.2.6}, one looks for $(a,b)\in\bC_*\times\bC$ such that
\begin{align}
aV_{\vz'}(a\ze+b)&=V_\vz(\ze)\label{eq4.5.1}\\
\sum_j\fr{\Gam_j}{\ze-\fr{z'_j-b}{a}}-a^2\ze-ab-ac&=\sum_j\fr{\Gam_j}{\ze-z_j}-\ze-c.\nonum
\end{align}
The equality between the two analytic parts already forces
\begin{equation}\label{eq4.5.2}
(a,b)=(1,0)\text{ \ or \ }(-1,-2c).
\end{equation}

\bigskip
In Section \ref{sec5.1}, through the simplest example of two vortices, we will illustrate that the bounds in Proposition \ref{pro4.2} and Theorem \ref{thm4.6} may or may not be attained. In this example, in addition, all the fixed equilibrium configurations come from the reduction, via Definition \ref{dfn4.4}, of the solutions of \eqref{eq4.2.1}, \eqref{eq4.2.2} or \eqref{eq4.2.3} by symmetries of the vortex sets and/or the given background flow only. Section \ref{sec5.2} will provide an example that a repeated solution of \eqref{eq4.2.1}, \eqref{eq4.2.2} or \eqref{eq4.2.3} exists, thus suggesting another factor in the reduction of the solutions to fixed equilibrium configurations via Definition \ref{dfn4.4}. All these results will be summarized in Tables \ref{tab5.1} and \ref{tab5.2}.

\bigskip
\subsection{Two Vortices}\label{sec5.1}

\bigskip
For two vortices (i.e. $n=2$), \eqref{eq4.2.3} reads
\begin{equation}\label{eq4.5.3}
\left\{\begin{matrix}
\begin{alignedat}{3}
&P_1(z_1,z_2)	&&=\Gam_1z_1+\Gam_2z_2+c(\Gam_1+\Gam_2)							&&=0\\
&P_2(z_1,z_2)	&&=\Gam_1{z_1}^2+\Gam_2{z_2}^2+c\,\Gam_1z_1+c\,\Gam_2z_2-\Gam_1\Gam_2	&&=0
\end{alignedat}\\
\Gam_1,\ \Gam_2\neq0
\end{matrix}\right.
\end{equation}
Computing by (the improved) Buchberger's algorithm (\cite[Algorithm GR\"OBNER-NEW2, p.232, \& Subalgorithm UPDATE, p.230]{BecWeKre}) with respect to the variable ordering $z_1>z_2$, we obtain a Gr\"obner basis
$$\cG=\cub{P_1,Q_3},\text{ \ where \ }Q_3:=-\tfr{\Gam_1+\Gam_2}{\Gam_1}{z_2}^2-\tfr{2c(\Gam_1+\Gam_2)}{\Gam_1}z_2+\Gam_1-\tfr{c^2(\Gam_1+\Gam_2)}{\Gam_1}\not\equiv0,$$
of $\cF=\cub{P_1,P_2}$. Thus, \eqref{eq4.5.3} is equivalent to $P_1=Q_3=0$, i.e.
\begin{equation}\label{eq4.5.4}
\left\{\begin{matrix}
\begin{alignedat}{2}
\Gam_1z_1+\Gam_2z_2+c(\Gam_1+\Gam_2)&=0\\
-\tfr{\Gam_1+\Gam_2}{\Gam_1}{z_2}^2-\tfr{2c(\Gam_1+\Gam_2)}{\Gam_1}z_2+\Gam_1-\tfr{c^2(\Gam_1+\Gam_2)}{\Gam_1}&=0
\end{alignedat}\\
\Gam_1,\ \Gam_2\neq0
\end{matrix}\right.
\end{equation}

\bigskip\noindent
\ul{Case 1:\quad$\Gam_1+\Gam_2=0$}, then $\Gam_1\neq\Gam_2$ because $\Gam_1,\Gam_2\neq0$, and then the second equation in \eqref{eq4.5.4} reads $\Gam_1=0$ which already has no solution. But this case is beyond the {\color{p}scope} of Proposition \ref{pro4.2} and Theorem \ref{thm4.6}.

\bigskip\noindent
\ul{Case 2:\quad$\Gam_1+\Gam_2\neq0$}, then the two solutions are
$$(z_1,z_2)=\rob{-c\mp\tfr{\Gam_2}{\surd(\Gam_1+\Gam_2)},-c\pm\tfr{\Gam_1}{\surd(\Gam_1+\Gam_2)}}=:\vz^\pm,$$
and the bound \smash{$\fr{(1+2-1)!}{(1-1)!}=2$} in Proposition \ref{pro4.2} is attained. Moreover, $\vz^\pm$ are admissible (see Remark (iv) following Definition \ref{dfn4.4}) since $\Gam_1+\Gam_2\neq0$, and are distinct since $\Gam_1,\Gam_2\neq0$.
\begin{enumerate}
\item[$\bullet$]	\ul{Case I:\quad$\Gam_1=\Gam_2$ (one species)}, then
\begin{align*}
2\pi\ii V_{\vz^-}(\ze)&=\fr{\Gam_1}{\ze+c-\sqrt{\tfr{\Gam_1}{2}}}+\fr{\Gam_1}{\ze+c+\sqrt{\tfr{\Gam_1}{2}}}-\ze-c\\
&=\fr{\Gam_1}{\ze+c+\sqrt{\tfr{\Gam_1}{2}}}+\fr{\Gam_1}{\ze+c-\sqrt{\tfr{\Gam_1}{2}}}-\ze-c=2\pi\ii V_{\vz^+}(\ze)
\end{align*}
so that \eqref{eq4.5.1} with $(a,b)=(1,0)$ in \eqref{eq4.5.2} is satisfied, and $\vz^\pm$ constitute only one fixed equilibrium configuration. In this case, the bound \smash{$\fr{(1+2-1)!}{(1-1)!2!}=1$} in Theorem \ref{thm4.6} is attained.
\item[$\bullet$]	\ul{Case II:\quad$\Gam_1\neq\Gam_2$ (two species)}, then 
\begin{align*}
2\pi\ii V_{\vz^-}(\ze)&=\fr{\Gam_1}{\ze+c-\tfr{\Gam_2}{\surd(\Gam_1+\Gam_2)}}+\fr{\Gam_2}{\ze+c+\tfr{\Gam_1}{\surd(\Gam_1+\Gam_2)}}-\ze-c\\
&\neq\fr{\Gam_1}{\ze+c+\tfr{\Gam_2}{\surd(\Gam_1+\Gam_2)}}+\fr{\Gam_2}{\ze+c-\tfr{\Gam_1}{\surd(\Gam_1+\Gam_2)}}-\ze-c=2\pi\ii V_{\vz^+}(\ze)
\end{align*}
so that \eqref{eq4.5.1} with $(a,b)=(1,0)$ in \eqref{eq4.5.2} is not satisfied, but
\begin{align*}
&-2\pi\ii V_{\vz^-}(-\ze-2c)=-\fr{\Gam_1}{-\ze-2c+c-\tfr{\Gam_2}{\surd(\Gam_1+\Gam_2)}}-\fr{\Gam_2}{-\ze-2c+c+\tfr{\Gam_1}{\surd(\Gam_1+\Gam_2)}}\\
&\qquad\qquad\qquad\qquad\qquad\qquad\qquad\qquad\quad+(-\ze-2c)+c\\
&\qquad\qquad\qquad\quad=\fr{\Gam_1}{\ze+c+\tfr{\Gam_2}{\surd(\Gam_1+\Gam_2)}}+\fr{\Gam_2}{\ze+c-\tfr{\Gam_1}{\surd(\Gam_1+\Gam_2)}}-\ze-c=2\pi\ii V_{\vz^+}(\ze)
\end{align*}
so that \eqref{eq4.5.1} with $(a,b)=(-1,-2c)$ in \eqref{eq4.5.2} is satisfied, hence $\vz^\pm$ still constitute only one fixed equilibrium configuration. In this case, the bound \smash{$\fr{(1+2-1)!}{(1-1)!}=2$} in Theorem \ref{thm4.6} is not attained.
\end{enumerate}

\begin{table}[h]
\centering
\begin{tabular}{|c|c|c|}
\hline
							&$\Gam_1+\Gam_2=0$	&$\Gam_1+\Gam_2\neq0$\\
\hline
$\Gam_1=\Gam_2$ (1 species)		&---					&$1$ / $2$ / $2$\\
\hline
$\Gam_1\neq\Gam_2$ (2 species)	&$0$ / $0$	 / $0$		&$1$ / $2$ / $2$\\
\hline
\end{tabular}
\caption{This table summarizes Section \ref{sec5.1}, i.e. the case $n=2$, $m=1$ and $W\equiv\text{constant}$ in \eqref{eq4.2.3}. `$a$ / $b$ / $c$' means that `the system \eqref{eq4.2.3} has $c$ solutions (counting multiplicity), $b$ of which are distinct and admissible, and that these solutions constitute $a$ fixed equilibrium configurations (in the sense of Definition \ref{dfn4.4})'. `---' means non-existence of the case.\label{tab5.1}}
\end{table}

\bigskip
\subsection{Repeated Solution}\label{sec5.2}

\bigskip
There is a case of three vortices in a quadrupole background flow (i.e. $n=3$ and \smash{$w(\ze)=-\fr{\ze}{2\pi\ii}$} in \eqref{eq4.2.1} or $m=1$ and $W\equiv0$ in \eqref{eq4.2.2} or \eqref{eq4.2.3}) where \eqref{eq4.2.1}, \eqref{eq4.2.2} or \eqref{eq4.2.3} has repeated solution. In such {\color{p}a} case, the bound in Theorem \ref{thm4.6} must not be attained. Consider \eqref{eq4.2.3} with the extra assumption that $\Gam_1=\Gam_2$ (at most two species):
\begin{equation}\label{eq4.5.6}
\left\{\begin{matrix}
\begin{alignedat}{3}
&P_1(z_1,z_2,z_3)	&&=\Gam_1z_1+\Gam_1z_2+\Gam_3z_3									&&=0\\
&P_2(z_1,z_2,z_3)	&&=\Gam_1{z_1}^2+\Gam_1{z_2}^2+\Gam_3{z_3}^2-{\Gam_1}^2-2\,\Gam_1\Gam_3	&&=0\\
&P_3(z_1,z_2,z_3)	&&=\Gam_1{z_1}^3+\Gam_1{z_2}^3+\Gam_3{z_3}^3-\Gam_1(\Gam_1+\Gam_3)z_1\\
&				&&\qquad\qquad\qquad-\Gam_1(\Gam_1+\Gam_3)z_2-2\,\Gam_1\Gam_3z_3			&&=0
\end{alignedat}\\
\Gam_1,\ \Gam_3\neq0
\end{matrix}\right.
\end{equation}
Computing by the aforesaid Buchberger's algorithm with respect to the variable ordering $z_1>z_2>z_3$, we obtain a Gr\"obner basis
$$\cG=\cub{P_1,Q_4,Q_5},$$
where
\begin{align*}
Q_4&:=-2\,\Gam_1{z_2}^2-2\,\Gam_3z_2z_3-\tfr{\Gam_3(\Gam_1+\Gam_3)}{\Gam_1}{z_3}^2+\Gam_1(\Gam_1+2\,\Gam_3)\not\equiv0\quad\text{and}\\
Q_5&:=-\tfr{\Gam_3(\Gam_1+\Gam_3)(2\,\Gam_1+\Gam_3)}{2\,{\Gam_1}^2}{z_3}^3+\tfr{\Gam_3(5\,\Gam_1+4\,\Gam_3)}{2}z_3\not\equiv0,
\end{align*}
of $\cF=\cub{P_1,P_2,P_3}$. Thus, \eqref{eq4.5.6} is equivalent to $P_1=Q_4=Q_5=0$, i.e.
\begin{equation}\label{eq4.5.7}
\left\{\begin{matrix}
\begin{alignedat}{2}
\Gam_1z_1+\Gam_1z_2+\Gam_3z_3&=0\\
-2\,\Gam_1{z_2}^2-2\,\Gam_3z_2z_3-\tfr{\Gam_3(\Gam_1+\Gam_3)}{\Gam_1}{z_3}^2+\Gam_1(\Gam_1+2\,\Gam_3)&=0\\
-\tfr{\Gam_3(\Gam_1+\Gam_3)(2\,\Gam_1+\Gam_3)}{2\,{\Gam_1}^2}{z_3}^3+\tfr{\Gam_3(5\,\Gam_1+4\,\Gam_3)}{2}z_3&=0
\end{alignedat}\\
\Gam_1,\ \Gam_3\neq0
\end{matrix}\right.
\end{equation}

\bigskip\noindent
\ul{Case 1:\quad$(\Gam_1+\Gam_3)(2\,\Gam_1+\Gam_3)=0$}, then $\Gam_1\neq\Gam_3$ because $\Gam_1,\Gam_3\neq0$, so that there are two species. And then writing $\Gam_3=-\al\Gam_1$, where $\al=1$ or $2$, \eqref{eq4.5.7} reads
$$\left\{\begin{matrix}
\begin{alignedat}{2}
\Gam_1z_1+\Gam_1z_2-\al\Gam_1z_3&=0\\
-2\,\Gam_1{z_2}^2+2\al\Gam_1z_2z_3+\al(1-\al)\Gam_1{z_3}^2+(1-2\al){\Gam_1}^2&=0\\
-\tfr{\al(5-4\al){\Gam_1}^2}{2}z_3&=0
\end{alignedat}\\
\Gam_1,\ \Gam_3\neq0
\end{matrix}\right.$$
The two solutions
$$(z_1,z_2,z_3)=\rob{\mp\sqrt{\tfr{(1-2\al)\Gam_1}{2}},\pm\sqrt{\tfr{(1-2\al)\Gam_1}{2}},0}=:\vz^\pm$$
are admissible and distinct, both because of $\Gam_1\neq0$ and \smash{$\al\neq\half$}. Moreover,
\begin{align*}
2\pi\ii V_{\vz^-}(\ze)&=\fr{\Gam_1}{\ze-\sqrt{\tfr{(1-2\al)\Gam_1}{2}}}+\fr{\Gam_1}{\ze+\sqrt{\tfr{(1-2\al)\Gam_1}{2}}}+\fr{-\al\Gam_1}{\ze}-\ze\\
&=\fr{\Gam_1}{\ze+\sqrt{\tfr{(1-2\al)\Gam_1}{2}}}+\fr{\Gam_1}{\ze-\sqrt{\tfr{(1-2\al)\Gam_1}{2}}}+\fr{-\al\Gam_1}{\ze}-\ze=2\pi\ii V_{\vz^+}(\ze)
\end{align*}
so that \eqref{eq4.5.1} with $(a,b)=(1,0)$ in \eqref{eq4.5.2} is satisfied, and they constitute only one fixed equilibrium configuration. But this case is beyond the {\color{p}scope} of Proposition \ref{pro4.2} and Theorem \ref{thm4.6}.

\bigskip\noindent
\ul{Case 2:\quad$(\Gam_1+\Gam_3)(2\,\Gam_1+\Gam_3)\neq0$}, then, first of all, the third equation of \eqref{eq4.5.7} has three solutions
$$z_3=0,\ \ \pm\Gam_1\sqrt{\tfr{5\,\Gam_1+4\,\Gam_3}{(\Gam_1+\Gam_3)(2\,\Gam_1+\Gam_3)}}.$$
Each leads to two $z_2$ via the second equation in \eqref{eq4.5.7} and then one $z_1$ via the first, resulting in six solutions of \eqref{eq4.5.7}, thus the bound \smash{$\fr{(1+3-1)!}{(1-1)!}=6$} in Proposition \ref{pro4.2} is attained. Now, one could observe that if
\begin{enumerate}
\item[$\bullet$]	\ul{$5\,\Gam_1+4\,\Gam_3=0$}, then $\Gam_1\neq\Gam_3$ because $\Gam_1,\Gam_3\neq0$, so that there are two species. Now, $z_3=0$ is actually a triple zero, and the six solutions of \eqref{eq4.5.7} are $(z_1,z_2,z_3)=$
$$\vz^\pm:=\rob{\mp\sqrt{\tfr{\Gam_1+2\,\Gam_3}{2}},\pm\sqrt{\tfr{\Gam_1+2\,\Gam_3}{2}},0}\quad\text{(each repeated thrice)}.$$
Note that \smash{$5\,\Gam_1+4\,\Gam_3=0\xRightarrow{\Gam_1,\Gam_3\neq0\ }\Gam_1+2\,\Gam_3\neq0$}, so $\vz^\pm$ are admissible and distinct. Moreover,
\begin{align*}
2\pi\ii V_{\vz^-}(\ze)&=\fr{\Gam_1}{\ze-\sqrt{\tfr{\Gam_1+2\,\Gam_3}{2}}}+\fr{\Gam_1}{\ze+\sqrt{\tfr{\Gam_1+2\,\Gam_3}{2}}}+\fr{-\fr{5}{4}\Gam_1}{\ze}-\ze\\
&=\fr{\Gam_1}{\ze+\sqrt{\tfr{\Gam_1+2\,\Gam_3}{2}}}+\fr{\Gam_1}{\ze-\sqrt{\tfr{\Gam_1+2\,\Gam_3}{2}}}+\fr{-\fr{5}{4}\Gam_1}{\ze}-\ze=2\pi\ii V_{\vz^+}(\ze)
\end{align*}
so that \eqref{eq4.5.1} with $(a,b)=(1,0)$ in \eqref{eq4.5.2} is satisfied, and they constitute only one fixed equilibrium configuration. Thus, the bound \smash{$\fr{(1+3-1)!}{(1-1)!2!}=3$} in Theorem \ref{thm4.6} is not attained.
\end{enumerate}

\begin{table}[h]
\centering
\begin{tabular}{|c|c|c|}
\hline
$\Gam_1=\Gam_2$				&$(\Gam_1+\Gam_3)(2\,\Gam_1+\Gam_3)=0$	&\begin{tabular}{@{}c@{}}$5\,\Gam_1+4\,\Gam_3=0$\\(necessarily,\\$(\Gam_1+\Gam_3)(2\,\Gam_1+\Gam_3)\neq0$)\end{tabular}\\
\hline
$\Gam_1=\Gam_3$ (1 species)		&---									&---\\
\hline
$\Gam_1\neq\Gam_3$ (2 species)	&$1$ / $2$ / $2$						&$1$ / $2$	 / $6$\\
\hline
\end{tabular}
\caption{This table summarizes Section \ref{sec5.2}, i.e. the case $n=3$, $m=1$ and $W\equiv0$ with $\Gam_1=\Gam_2$ in \eqref{eq4.2.3}. `$a$ / $b$ / $c$' means that `the system \eqref{eq4.2.3} has $c$ solutions (counting multiplicity), $b$ of which are distinct and admissible, and that these solutions constitute $a$ fixed equilibrium configurations (in the sense of Definition \ref{dfn4.4})'. `---' means non-existence of the case.\label{tab5.2}}
\end{table}



\bigskip
\begin{thebibliography}{99}
{\color{g}
\bibitem{Ac}
D. J. ACHESON, {\it Elementary Fluid Dynamics}, Clarendon Press, Oxford, 1990.

\bibitem{BecWeKre}
T. BECKER, V. WEISPFENNING and H. KREDEL, {\it Gr\"obner Bases: A Computational Approach to Commutative Algebra}, Springer-Verlag, New York, 1993.

\bibitem{Ber}
D. N. BERNSHTEIN. The number of roots of a system of equations, {\it Functional Anal. Appl.} 9:183-185 (1975). Translated from {\it Funkts. Anal. Prilozh.} 9:1-4 (1975).

\bibitem{CaKad}
L. J. CAMPBELL and J. B. KADTKE. Stationary configurations of point vortices and other logarithmic objects in two dimensions, {\it Phys. Rev. Lett.} 58:670-673 (1987).

\bibitem{Cl}
P. A. CLARKSON. Vortices and polynomials, {\it Stud. Appl. Math.} 123:37-62 (2009).

\bibitem{Cox}
D. A. COX. Solving equations via algebras, In {\it Solving Polynomial Equations: Foundations, Algorithms, and Applications}, edited by A. Dickenstein, I. Z. Emiris et al., 63-123, Springer, Berlin, 2005.

\bibitem{CoxLitOS2}
D. A. COX, J. B. LITTLE and D. O'SHEA, {\it Using Algebraic Geometry}, Springer, New York, 2005.

\bibitem{Ha}
M. HAMPTON. Finiteness of kite relative equilibria in the five-vortex and five-body problems, {\it Qual. Theory Dyn. Syst.} 8:349-356 (2009).

\bibitem{HaMo1}
M. HAMPTON and R. MOECKEL. Finiteness of relative equilibria of the four-body problem, {\it Invent. Math.} 163:289-312 (2006).

\bibitem{HaMo2}
M. HAMPTON and R. MOECKEL. Finiteness of stationary configurations of the four-vortex problem, {\it Trans. Amer. Math. Soc.} 361:1317-1332 (2008).

\bibitem{KadCa}
J. B. KADTKE and L. J. CAMPBELL. Method for finding stationary states of point vortices, {\it Phys. Rev. A} 36:4360-4370 (1987).

\bibitem{Kh}
A. G. KHOVANSKII. Newton polyhedra and the genus of complete intersections, {\it Functional Anal. Appl.} 12:38-46 (1978). Translated from {\it Funkts. Anal. Prilozh.} 12:51-61 (1978).

\bibitem{Kus}
A. G. KUSHNIRENKO. Newton polytopes and the Bezout theorem, {\it Functional Anal. Appl.} 10:233-235 (1976). Translated from {\it Funkts. Anal. Prilozh.} 10:82-83 (1976).

\bibitem{LiWa}
T. Y. LI and X.-S. WANG. The BKK root count in $\vC^n$, {\it Math. Comp.} 65:1477-1484 (1996).

\bibitem{MeAr}
V. V. MELESHKO and H. AREF. A bibliography of vortex dynamics 1858--1956, {\it Adv. Appl. Mech.} 41: 197-292 (2007).

\bibitem{ONe1}
K. A. O'NEIL. Stationary configurations of point vortices, {\it Trans. Amer. Math. Soc.} 302:383-425 (1987).

\bibitem{ONe2}
K. A. O'NEIL. Minimal polynomial systems for point vortex equilibria, {\it Phys. D} 219:69-79 (2006).
}
\end{thebibliography}
\end{document}